\documentclass[11pt]{article}

\usepackage{amsmath}
\usepackage{amssymb}
\usepackage{amscd}
\usepackage{latexsym}
\usepackage{theorem}
\usepackage{setspace}
\usepackage{epsfig}
\usepackage{psfrag}
\usepackage{amsfonts}
\usepackage{enumerate}
\usepackage{euscript}

\newtheorem{Def}{Definition}[section]
\newtheorem{theorem}{Theorem}[section]
\newtheorem{lemma}[theorem]{Lemma}
\newtheorem{proposition}[theorem]{Proposition}
\newtheorem{corollary}[theorem]{Corollary}

{\theorembodyfont{\rmfamily} \theoremstyle{plain}

\newtheorem{example}[theorem]{Example}

\newtheorem{remark}[theorem]{Remark}

}

\numberwithin{equation}{section} \numberwithin{figure}{section}

\newcommand{\CP}{\mathbb{CP}}
\newcommand{\C}{\mathbb{C}}

\newcommand{\RP}{\mathbb{RP}}

\newcommand{\md}[2]{\ensuremath{\overline{M}_{#1}( #2 ,d)}}

\begin{document}

\noindent {\LARGE \bf Real Gromov-Witten invariants on the \\
  moduli space of genus 0 stable maps to \\  \indent a
smooth rational projective space }
\begin{center}
\Large{\emph{Dedicated to the originator Gang Tian}}
\end{center}

\begin{center}
{\bf Seongchun Kwon }\footnote{
\newline \mbox{~~~~}{\it MSC 2000 Subject
Classification}: Primary: 14C17 \hspace{0.1in} Secondary: 14C25
\newline \mbox{~~~~}{\it Keywords}: Gromov-Witten invariant,
enumerative invariant, transversality, intersection multiplicity,
real structure}
\end{center}

{\small
\begin{quote}
\noindent {\em Abstract.} We characterize transversality,
non-transversality properties on the moduli space of genus 0 stable
maps to a rational projective surface.
 If a target space is equipped with a real structure, i.e, anti-holomorphic involution, then
the results have real enumerative applications. Firstly, we can
define a real version of Gromov-Witten invariants. Secondly, we can
prove the invariance of Welschinger's invariant in algebraic
geometric category.
\end{quote}}

\bigskip

\section{Introduction} \label{s;intro}

Let $\overline{M}_{k}(X,\beta)$ be the moduli space of stable maps
from a $k$-pointed arithmetic genus 0 curve to $X$, representing a
2nd homology class $\beta$. Let $[\Upsilon_{1}]$, \ldots,
$[\Upsilon_{k}]$ be Poincar\'{e} duals to the homology classes
represented by $\Upsilon_{1}$, \ldots, $\Upsilon_{k}$, where
$\Upsilon_{1}$, \ldots, $\Upsilon_{k}$ are pure dimensional
varieties in the target space $X$. The Gromov-Witten invariant on
$\overline{M}_{k}(X,\beta)$ is defined as:
\begin{center}
$\displaystyle I_{\beta}([\Upsilon_{1}], \ldots, [\Upsilon_{k}]) :=
  \int_{\overline{M}_{k}(X,\beta)} ev_{1}^{*}([\Upsilon_{1}]) \cup
\ldots \cup ev_{k}^{*}([\Upsilon_{k}])$,
\end{center}
where $ev_{i}$ is an $i$-th evaluation map. \\
The Gromov-Witten invariant $I_{\beta}([\Upsilon_{1}], \ldots,
[\Upsilon_{k}])$ may be non-trivial only when $\sum \mbox{codim}
(\Upsilon_{i}) = \mbox{dim} \overline{M}_{k}(X,\beta)$. We say that
\emph{the Gromov-Witten invariant has an enumerative meaning} if
$I_{\beta}([\Upsilon_{1}], \ldots, [\Upsilon_{k}])$ equals to the
actual number of points in $ev_{1}^{-1}(\Gamma_{1}) \cap \ldots \cap
ev_{k}^{-1}(\Gamma_{k})$, where $\Gamma_{1}$, \ldots, $\Gamma_{k}$
are any pure dimensional varieties in a general position such that
$[\Gamma_{i}] = [\Upsilon_{i}]$, $i=1, \ldots, k$. So, the
Gromov-Witten invariant counts the number of stable maps whose
$i$-th marked point maps into $\Gamma_{i}$ if it has an enumerative
implication. Note that the number of intersection points in
$ev_{1}^{-1}(\Gamma_{1}) \cap \ldots \cap ev_{k}^{-1}(\Gamma_{k})$
doesn't vary depending on the general choices of the cycle's
representatives $\Gamma_{i}$. And $ev_{i}^{-1}(\Gamma_{i})$, $i=1,
\ldots, k$, meet transversally for the general choices of the
cycle's representatives $\Gamma_{i}$. As it is well-known, the
Gromov-Witten invariant has an enumerative implication if the target
space is a homogeneous variety. See \cite{fap}. See \cite{gap} for
an enumerative application of the Gromow-Witten invariant of
blow-ups of $\CP^{2}$ at the finite number of points.
\par Let $X$ be a rational projective surface. That is, $X$ is
deformation equivalent to either $\CP^2$, $\CP^1 \times \CP^1$ or
$\CP^2$ blown-up at the finite number of points $B := \{b_{1},
\ldots, b_{r}\}$ which we will denote by $r \CP^2$. The aim of this
paper is to study the intersection theoretic properties on
$\overline{M}_{k}(X,\beta)$ and real enumerative applications when
$X$ is equipped with a real structure induced by a complex
conjugation map on the complex projective space. That is,
\begin{itemize}
\item $\CP^2 \rightarrow \CP^2$, $[a:b:c] \mapsto [\overline{a}:
\overline{b}: \overline{c}]$
\item $\CP^1 \times \CP^1 \rightarrow \CP^1 \times \CP^1$,
$([a:b],[c:d]) \mapsto ([\overline{a}:\overline{b}], [\overline{c},
\overline{d}])$
\end{itemize}

Let $r \CP^2$ come from blown-ups of $\CP^2$ at $B := \{b_{1},
\ldots, b_{r} \}$, where $B$ is preserved by a complex conjugation
map on $\CP^2$. Then, $r \CP^2$ has an obvious real structure which
is induced by a complex conjugation involution on $\CP^2$.
\par \emph{A real projective variety} is a projective variety defined over $\C$,
 having a real structure, i.e., an anti-holomorphic involution. \emph{A real
 part} of the real projective variety is the locus which is fixed by an
 anti-holomorphic involution.
Let $X$ be a real projective rational surface which is described
above. Then, $\overline{M}_{k}(X,\beta)$ is a real projective
variety. Based on the results in ~\cite{kwon}, we can study the real
enumerative problems on the real part
$\overline{M}_{k}(X,\beta)^{re}$.
\par Let $\Gamma_{1}$, \ldots, $\Gamma_{k}$ be any pure dimensional
\emph{real} projective varieties in the real rational surface $X$ in
a general position such that $[\Gamma_{i}] = [\Upsilon_{i}]$, $i=1,
\ldots, k$. Assume that the Gromov-Witten invariant
$I_{\beta}([\Upsilon_{1}], \ldots, [\Upsilon_{k}])$ is non-trivial.
Then, the number of intersection points $ev_{1}^{-1}(\Gamma_{1})
\cap \ldots \cap ev_{k}^{-1}(\Gamma_{k}) \cap
\overline{M}_{k}(X,\beta)^{re}$ varies depending on the actual
choice of $\Gamma_{1}, \ldots, \Gamma_{k}$. However, if each point
in $ev_{1}^{-1}(\Gamma_{1}) \cap \ldots \cap ev_{k}^{-1}(\Gamma_{k})
\cap \overline{M}_{k}(X,\beta)^{re}$ has an intersection
multiplicity one, then the number of points doesn't vary by little
perturbations of cycle's representatives $\Gamma_{1}, \ldots,
\Gamma_{k}$. The changes of the number of intersection points in
$ev_{1}^{-1}(\Gamma_{1}) \cap \ldots \cap ev_{k}^{-1}(\Gamma_{k})
\cap \overline{M}_{k}(X,\beta)^{re}$ happen only after some of the
intersection points have intersection multiplicities greater than
one. So, it is important to study  transversality and
non-transversality properties.

\begin{Def} Let $X$ be deformation equivalent to $\CP^{2}$, $\CP^1
\times \CP^1$ or $\CP^2$ blown-up at finite number of points.\\
(a) A stable map $f: C \rightarrow X$ is called a \emph{cuspidal
stable map} if $C$ is isomorphic to $\CP^1$, and its image $f(C)$
contains only node singularities and a unique cuspidal
singularity.\\
(b) By an \emph{equi-singular locus}, we mean the set of stable maps
in $\overline{M}_{k}(X,\beta)$ which have the same type of
singularities on the image curves if the domain curve is
irreducible. If the domain curve is reducible, then an
\emph{equi-singular locus} means the set of stable maps having the
same number of irreducible components in the domain curve.
\end{Def}

\par Theorem \ref{t;main} and Theorem
\ref{t;mainblow} are the main results of this paper. Both Theorems
are important for real enumerative applications. However, in both
Theorems, we don't need to assume that $X$ has a real structure.

\begin{theorem}\label{t;main}
Let $X$ be deformation equivalent to either $\CP^2$ or $\CP^1 \times
\CP^1$. Let $\Gamma_{i}$, $i=1, \ldots,k$, be pure dimensional
subvarieties in $X$. Let $p$ be a point in $ev_{1}^{-1}(\Gamma_{1})
\cap \ldots \cap ev_{k}^{-1}(\Gamma_{k})$, which represents a
cuspidal stable map. Assume that $\displaystyle \sum_{i=1}^{k}
\hspace{1mm}
\mbox{codim} \Gamma_{i}= \mbox{dim} \overline{M}_{k}(X, \beta)$. \\
(i) Assume $k= -\int_{\beta} \omega_{X}-1$. Then, the intersection
multiplicity at $p$ is 2. The cuspidal stable maps locus is a unique
equi-singular locus having a codimension $\leq 1$ on which a
transversality always
fails. \\
(ii) Assume $k > -\int_{\beta} \omega_{X}-1$. Suppose that the image
of the stable map represented by $p$ meets all non-trivial Chow
cycle's representatives
 $\Gamma_{i}$ transversally. Then, the intersection multiplicity at $p$ is 2.
Transversality uniformly fails along the cuspidal stable maps'
locus.
\end{theorem}

\par An easy consequence of Theorem ~\ref{t;main}
 is the number of points
in $ev_{1}^{-1}(\Gamma_{1}) \cap \ldots \cap
ev_{k}^{-1}(\Gamma_{k})$ is always strictly less than the
Gromov-Witten invariant $I_{\beta}([\Gamma_{1}], \ldots,
[\Gamma_{k}])$ if any one of the points in $ev_{1}^{-1}(\Gamma_{1})
\cap \ldots \cap ev_{k}^{-1}(\Gamma_{k})$ is a cuspidal stable map.

\par A variety $X$ is called a \emph{convex variety} if $H^{1}(C, f^{*}
TX)$ vanishes for all arithmetic genus 0 curves $C$. Let $r \CP^2$
be deformation equivalent to $\CP^2$ blown up at $r$ points
$\{b_{1}, \ldots, b_{r} \}$. Then, $r \CP^2$ is a non-convex
variety. Although we need to consider a virtual fundamental cycle,
the Gromov-Witten invariant has an enumerative meaning. See
\cite{gap}. That is, the intersection theory on the reducible domain
curve's locus is not important for the enumerative application
purpose in complex category. For a real enumerative application, we
consider only a \emph{non-exceptional divisor class} $d \cdot$[line]
in $r \CP^2$ in this paper.

\begin{theorem}\label{t;mainblow}
Let $p_{i}$, $i=1, \ldots, 3d-1$, be points in general position in
$r \CP^2 \setminus (\CP^{1}_{1} \cup \ldots \cup \CP^{1}_{r} )$,
where
 $\CP^{1}_{j}$ is an exceptional divisor, $j=1, \ldots, r$.
Let $p$ be a point in $ev_{1}^{-1}(p_{1}) \cap \ldots \cap
ev_{3d-1}^{-1}(p_{3d-1}) \subset M_{3d-1}( r \CP^2, d \cdot
[line])$, which represents a cuspidal stable map. Then, the
intersection multiplicity at $p$ is 2. A cuspidal stable maps locus
is a unique equi-singular locus in $M_{3d-1}( r \CP^2, d \cdot
[line])$ which has a codimension $\leq 1$ and on which a
transversality always fails.
\end{theorem}

\par When the target space $X$ is equipped with a real structure,
Theorem ~\ref{t;main}, Theorem ~\ref{t;mainblow} play key roles in
defining real Gromov-Witten invariants which are local invariants on
$X^{re} \times \ldots \times X^{re}: = (X^{k})^{re}$, in
$I_{\beta}([\mbox{point}], \ldots, [\mbox{point}])$ case. Other
natural application is the proof of the invariance of the
Welschinger's invariant in algebraic geometric category.

\par One of the possible application problems suggested by Gang Tian
is the following problem. I invite the challenging readers to
attempt to resolve the following real enumerative problem:

\vspace{2mm}

\par Prove or disprove: There are 11 real configuration points in
$\CP^2$ such that all degree 4 rational nodal curves passing through
those configuration points are real curves.

\vspace{2mm}

\par The result with 564 real configuration points in \cite{mik} is the currently best toward this
problem. See Theorem 3.6 in ~\cite{sot} for the degree 3 case.\\

\par  The paper is organized as follows:\\
In Sec.\ref{s;ns}, we show that the nodal Severi variety is embedded
into $\md{0}{\CP^2}$ as an open locus on which the intersection
theory is established. Then, by using the classical results, we
classify the degeneration properties of the nodal stable maps in
$\md{0}{\CP^2}$. In Sec.\ref{ss;gromov}, we go through the local
calculations of the differential of the $ev$ map on the loci whose
codimensions are less than or equal to one in $\md{3d-1}{\CP^2}$.
The index of the $ev$ map is identical to the intersection
multiplicity. In Sec.\ref{ss;point}, we characterize walls and
chambers in $\md{3d-1}{\CP^2}$ and the Chow 0-cycles parameter
space. And then, we define the real version of Gromov-Witten
invariants in $I_{d}$([point], \ldots, [point]) case. In Sec.
\ref{ss;def}, we work on the $I_{d}$([point], \ldots,
[point],[line], \ldots, [line]) case. Different from the results in
Sec. \ref{ss;point}, walls and chambers are not characterized by the
geometric properties of the stable maps. In general case, we need to
add the tangency condition in our consideration. In
sec.\ref{s;rcp2}, we work on $r \CP^2$ with the divisor class of $d
\cdot $[line] and $\CP^1 \times \CP^1$ target space case.          \\

\par Throughout the paper, we will assume the degree of stable maps is greater
than or equal to 3 if the target space is $\CP^2$ or $r \CP^2$. But
tangent space splitting calculations done in Section~\ref{ss;gromov}
hold for all degree $d$.
 We will use a notation $d$ when we consider $d \cdot \mbox[line]$
in most cases.

\vspace{2mm}

\section{Transversality properties on
$\md{k}{\CP^{2}}$ and their real enumerative
implications}\label{sec;cp2}

\subsection{The classical rational nodal Severi variety v.s. the moduli space
of stable maps $\md{0}{\CP^{2}}$} \label{s;ns}
\par \emph{The geometric genus of a plane nodal curve} is the sum of
an arithmetic genus on each component after we desingularize all
nodes.

\par \emph{A rational nodal
Severi variety $\mathcal{NS}^{d}$ of degree $d$} is the set of
reduced, irreducible plane curves of degree $d$, having only nodes
as singularities, and having geometric genus zero. The rational
nodal Severi variety $\mathcal{NS}^{d}$ is a $3d-1$ dimensional
(quasi-)projective variety.

\begin{proposition}\label{p;dense}
The rational nodal Severi variety $\mathcal{NS}^{d}$ of degree $d$
is embedded into the moduli space $\md{0}{\CP^2}$ of stable maps as
an open sublocus.
\end{proposition}

Proof. Let $\textbf{c}$ be a point in $\mathcal{NS}^{d}$. Let $C$ be
a reduced, irreducible, rational nodal curve in $\CP^2$ which is
represented by $\textbf{c}$. Consider the morphism $F$ from
$\mathcal{NS}^{d}$ to $\md{0}{\CP^{2}}$ defined by $\textbf{c}
\mapsto [(f, \widetilde{C})]$, where $f$ is a normalization of $C$.
The morphism is well-defined because the normalization is unique up
to isomorphism and the isomorphism is exactly the equivalence
relation of stable maps.
\par Let $T_{\textbf{c}} \mathcal{NS}^{d}$ be a tangent space at
$\textbf{c} \in \mathcal{NS}^{d}$. Then, the Kodaira-Spencer map
$\mathcal{K}_{\textbf{c}}: T_{\textbf{c}} \mathcal{NS}^{d}
\rightarrow H^{0}(\widetilde{C}, N_{\textbf{c}})$ is onto, where
$N_{\textbf{c}}$ is Coker($df: T \widetilde{C} \rightarrow f^{*}
T\CP^2)$. See \cite[p110]{ham}. $\mathcal{K}_{\textbf{c}}$ is
isomorphic because the degree of the line bundle $N_{\textbf{c}}$ is
$3d-2$. The standard $K$-group calculations from the tangent
obstruction long exact sequence show the following isomorphism of
the tangent space at $[(f, \widetilde{C})] \in \md{0}{\CP^2}$ :
\begin{gather} \label{e;kod}  T_{[(f, \widetilde{C})]} \md{0}{\CP^2} \cong
\mbox{Ext}^{1}(f^{*} \Omega^{1}_{\CP^{2}} \rightarrow
\Omega^{1}_{\widetilde{C}}, \mathcal{O}_{\widetilde{C}}) \cong
H^{0}(\widetilde{C}, N_{\textbf{c}})
\end{gather}

Thus, the Kodaira-Spencer map with (\ref{e;kod}) shows that the
differential $dF_{\textbf{c}}$ at $\textbf{c}$ is an isomorphism.
The Proposition follows because $\mathcal{NS}^{d}$ is a
quasi-projective
variety. \hfill q.e.d.\\

In fact, the rational nodal Severi variety $\mathcal{NS}^{d}$ is
embedded into the fine moduli locus in $\md{0}{\CP^2}$. The reason
is there is a universal curve over $\mathcal{NS}^{d}$ such that the
fiber over each point in $\mathcal{NS}^{d}$ is the normalization of
the corresponding plane curve at the node singularities and the
canonical morphism from the universal curve to $\CP^2$. See p6 in
\cite{dah}.

\begin{corollary}\label{c;dense}
Let $\mathcal{NL} \subset \md{k}{\CP^2}$ be the locus of pointed
stable maps $[(f, \CP^1, a_{1}, \ldots, a_{k})]$ which satisfies the
following:
\begin{itemize}
\item The image curve $f(\CP^1)$ is reduced, irreducible, rational nodal curves
\item $f(a_{1}), \ldots, f(a_{k})$ don't meet the node singularities in $f(\CP^1)$.
\end{itemize}
Then, $\mathcal{NL} $ is an open sublocus.
\end{corollary}

Proof. Note that the condition in the second item is an open
condition. Since the forgetful map is continuous, Proposition
\ref{p;dense} implies that $\mathcal{NL} $ is an open sublocus.
\hfill q.e.d.\\

\begin{remark}
As it is well-known, the Gromov-Witten invariant counts irreducible
rational nodal curves passing through $3d-1$ points in general
position in $\CP^{2}$. The proof of Proposition \ref{p;dense} shows
the exact correspondences between the irreducible rational nodal
curves and the stable maps passing through the same set of $3d-1$
points. Therefore, the locus in Corollary \ref{c;dense} is the locus
where the intersection theory is established.
\end{remark}

\begin{figure}
\includegraphics{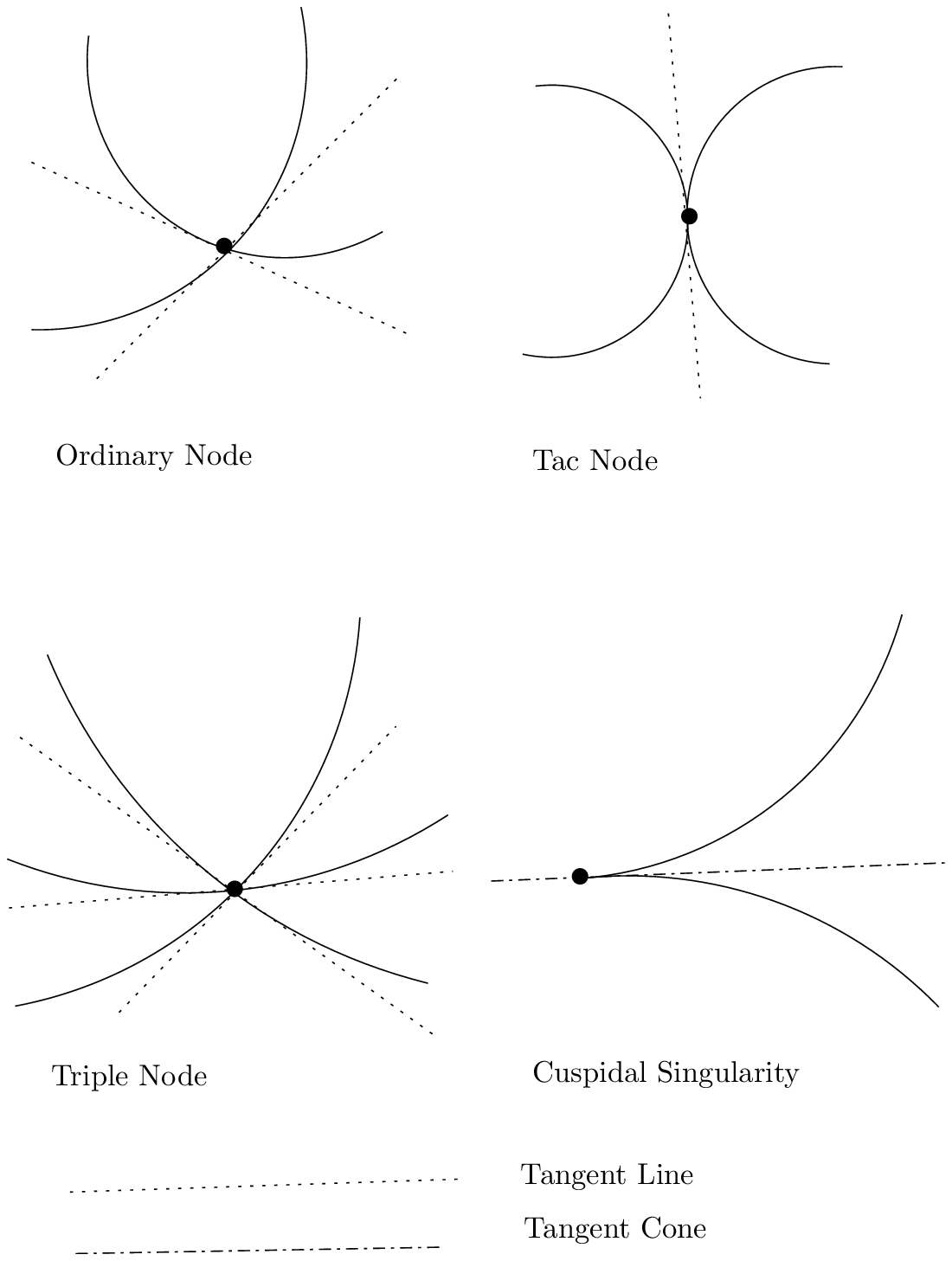}
\caption{Type of singularities} \label{p;nodetype}
\end{figure}

\begin{Def}
\emph{A nodal stable map} is a stable map $f:C \rightarrow \CP^2$
such that $C$ is isomorphic to $\CP^1$ and its image $f(C)$ contains
only (ordinary) node singularities.
\end{Def}

\begin{Def}
\emph{A tacnode stable map} is a stable map $f: C \rightarrow \CP^2$
such that $C$ is isomorphic to $\CP^1$ and its image $f(C)$ contains
only (ordinary) node singularities and
a unique tacnode.\\
\indent \emph{A triple node stable map} is a stable map $f: C
\rightarrow \CP^2$ such that $C$ is isomorphic to $\CP^1$ and its
image $f(C)$ contains only (ordinary) node singularities and a
unique triple node.
\end{Def}

\par The following Theorem characterizes the codimension one
equi-singular locus. One can easily see that the same
characterization holds in $\md{k}{\CP^2}$.
\par In the proof of Theorem \ref{t;stdeg}, the cuspidal curve locus in the partially compactified nodal
Severi variety $\mathcal{NS}^{d}$ is the locus where each curve has
exactly one cuspidal singularity and all other singularities are
node singularities.

\begin{theorem}\label{t;stdeg}
Let $\mathcal{NL}$ denote the nodal stable maps locus in
$\md{0}{\CP^2}$. The codimension one equi-singular loci in
$\md{0}{\CP^{2}} \setminus \mathcal{NL}$, are the cuspidal stable
maps locus, the tacnode stable maps locus, the triplenode stable
maps locus, the locus consisting of the stable maps whose domain
curves consist of two irreducible components.
\end{theorem}

Proof. The last case is well-known.
\par ~\cite[(1.4) Theorem]{dah2}
shows that the partially compactified rational nodal Severi variety
has the same classification of the codimension one equi-singular
loci with the classification of the codimension one equi-singular
loci in $\md{0}{\CP^2}$ in this Theorem. There is a universal curve
on the normalized partially complactified rational nodal Severi
variety, whose fiber over each point is the normalization of the
corresponding plane curve at the assigned singularities. And there
is a canonical morphism from the universal curve to $\CP^2$. See p6
in \cite{dah}. Thus, there is a canonical embedding from the
normalized partially compactified rational nodal Severi variety to
the fine moduli locus in $\md{0}{\CP^2}$. This canonical embedding
relates a plane curve represented by a point in the partially
compactified rational nodal Severi variety with the normalization of
the plane curve in $\md{0}{\CP^2}$. The normalization on the
partially compactified rational nodal Severi variety doesn't change
the codimension properties in the partially compactified nodal
Severi variety. Thus, the Theorem follows.  \hfill q.e.d.\\

\subsection{Transversality property on $\md{k}{\CP^{2}}$ and
its real enumerative implications}\label{s;gwcp2} \vspace{1mm}

\subsubsection{Transversality properties for the Gromov-Witten invariant
$I_{d}$([point], \ldots, [point]) case } \label{ss;gromov}

\vspace{1mm}

\par The tangent space at $[(f,\CP^1, a_{1}, \ldots, a_{k})]$
 in $\md{k}{\CP^2}$ is the hyperext group $Ext^{1}(f^{*}\Omega^{1}_{\CP^2} \rightarrow
\Omega^{1}_{\CP^{1}}(a_{1}+ \ldots + a_{k}),
\mathcal{O}_{\CP^{1}})$. We calculate the hyperext group in terms of
the ordinary sheaf cohomology group.

\vspace{2mm}

\begin{lemma}\label{l;tan}
Let $[(f,\CP^1, a_{1}, \ldots, a_{k})]$ be a point in
$\md{k}{\CP^2}$. Then, the tangent space at  $[(f,\CP^1, a_{1},
\ldots, a_{k})]$ is:
\begin{gather*}
 T \md{k}{\CP^2} \mid_{[(f, \CP^1, a_{1}, \ldots, a_{k})]} \hspace{1mm}
\cong \hspace{1mm} H^{0}(\CP^1, N_{f}) \oplus T_{a_{1}} \CP^1 \oplus
\ldots \oplus T_{a_{k}} \CP^{1}
\end{gather*}
where $N_{f}$ is Coker$(df: T \CP^1 \rightarrow f^{*} T \CP^2)$.
\end{lemma}
\vspace{2mm}

Proof. \hfill The \hfill long \hfill exact \hfill sequence
\hfill associated \hfill to \hfill the \hfill hyperext \hfill group\\
$Ext^{1}(f^{*}\Omega^{1}_{\CP^2} \rightarrow
\Omega^{1}_{\CP^{1}}(a_{1}+ \ldots + a_{k}), \mathcal{O}_{\CP^{1}})$
is:
\begin{multline} \label{e;vir}
0 \hspace{1mm} \rightarrow \hspace{1mm} Hom(\Omega_{\CP^1}^{1}(a_{1}
+ \ldots + a_{k}), \mathcal{O}_{\CP^1}) \hspace{1mm} \rightarrow
\hspace{1mm} H^{0}(
\CP^1, f^{*}T \CP^{2}) \rightarrow \\
\rightarrow Ext^{1}(f^{*}\Omega^{1}_{\CP^{2}} \rightarrow
\Omega^{1}_{\CP^{1}}(a_{1} + \ldots + a_{k}),
\mathcal{O}_{\CP^{1}})\rightarrow \\
 \rightarrow
Ext^{1}(\Omega^{1}_{\CP^{1}}(a_{1} + \ldots
+a_{k}),\mathcal{O}_{\CP^{1}}) \rightarrow 0
\end{multline}

The last term $Ext^{1}(\Omega^{1}_{\CP^{1}}(a_{1} + \ldots
+a_{k}),\mathcal{O}_{\CP^{1}}) $ of the deformation of a pointed
curve is isomorphic to $H^{1}(\CP^{1}, T \CP^{1}(-a_{1}- \ldots
-a_{k}))$. It comes from an exact sequence of the local to global
spectral sequence
\begin{multline}\label{e;dimu}
0 \hspace{1mm} \rightarrow \hspace{1mm} H^{1}(\CP^1,
\underline{Ext}^{0}_{\CP^{1}}(\Omega^{1}_{\CP^{1}}( a_{1} + \ldots
+ a_{k}), \mathcal{O}_{\CP^1})) \hspace{1mm} \rightarrow  \\
\rightarrow \hspace{1mm} Ext^{1}(\Omega^{1}_{\CP^1}(a_{1}+ \ldots +
a_{k}), \mathcal{O}_{\CP^{1}}) \hspace{1mm} \rightarrow \\
\rightarrow \hspace{1mm} H^{0}(\CP^1,
\underline{Ext}^{1}_{\CP^{1}}(\Omega^{1}_{\CP^{1}}(a_{1} + \ldots +
a_{k}),\mathcal{O}_{\CP^{1}})) \hspace{1mm} \rightarrow \hspace{1mm}
0
\end{multline}
by noting

\begin{itemize}
\item  $H^{0}(\hspace{1mm} \CP^1, \hspace{1mm}
\underline{Ext}^{1}_{\CP^{1}}( \hspace{1mm} \Omega^{1}_{\CP^{1}}(
\hspace{1mm} a_{1} + \ldots + a_{k} \hspace{1mm} ),\hspace{1mm}
\mathcal{O}_{\CP^{1}}\hspace{1mm} )\hspace{1mm} )$
\hfill is \hfill isomorphic \hfill to \\
$H^{0}(\CP^1,\underline{Ext}^{1}_{\CP^{1}}(\Omega^{1}_{\CP^{1}},\mathcal{O}_{\CP^{1}}))$
because the first order of smoothing is supported by nodes and
independent of marked points. Thus, it vanishes.
\item $\underline{Ext}^{0}_{\CP^1}( \Omega_{\CP^{1}}^{1}(a_{1} + \ldots +
a_{k}), \mathcal{O}_{\CP^1})$ is isomorphic to
$Hom(\Omega_{\CP^{1}}^{1}(a_{1} + \ldots +
a_{k}),\mathcal{O}_{\CP^{1}})$.
 Hence,  $H^{1}(\CP^1,
\underline{Ext}^{0}_{\CP^{1}}( \Omega^{1}_{\CP^{1}}( a_{1} + \ldots
+ a_{k} ), \mathcal{O}_{\CP^1} ) )$
 is isomorphic to
$H^{1}(\CP^{1},T \CP^{1}(-a_{1} - \ldots - a_{k}))$.
\end{itemize}

\par Thus, we get a splitting of a tangent space
$Ext^{1}(f^{*}\Omega^{1}_{\CP^{2}} \rightarrow
\Omega^{1}_{\CP^{1}}(a_{1} + \ldots + a_{k}),
\mathcal{O}_{\CP^{1}})$ which is naturally isomorphic to $ \ominus
Hom(\Omega^{1}_{\CP^{1}}(a_{1} + \ldots + a_{k}),
\mathcal{O}_{\CP^{1}}) \oplus H^{0}(\CP^{1}, f^{*}(T \CP^{2}))
\oplus H^{1}(\CP^{1}, T \CP^{1}(-a_{1} - \ldots - a_{k}))$.

\vspace{2mm}

\par The result follows from $K$-group calculations associated to a long exact
sequence induced from the following two short exact sequences of
sheaves:

\begin{equation} 0 \rightarrow T \CP^{1}(-a_{1} - \ldots - a_{k})
\rightarrow T \CP^{1} \rightarrow T_{a_{1}} \CP^{1} \oplus \ldots
\oplus T_{a_{k}} \CP^{1} \rightarrow 0 \label{e;sh}
\end{equation}
\begin{equation}
0 \rightarrow T \CP^{1} \rightarrow f^{*} T \CP^{2} \rightarrow
N_{f} \rightarrow 0 \label{e;nor1}
\end{equation}

More precisely,\\

\begin{align}
 &Hom(\Omega^{1}_{\CP^1}(a_{1} + \ldots + a_{k}),
\mathcal{O}_{\CP^1}) \notag\\
\cong &  H^{0}(\CP^1, T \CP^{1}( -a_{1} - \ldots - a_{k}))\notag\\
\cong & H^{0}( \CP^1, T \CP^1) \ominus ( T_{a_{1}}\CP^{1} \oplus
\ldots \oplus T_{a_{k}} \CP^1) \notag\\ & \oplus H^{1}( \CP^1, T
\CP^{1}( -a_{1} - \ldots - a_{k})) \hspace{2mm} \text{ by
(\ref{e;sh})}\notag
\end{align}

$H^{0}(\CP^1, f^{*}( T\CP^{2})) \cong H^{0}( \CP^1, T \CP^1) \oplus
H^{0}( \CP^1, N_{f})$ by (\ref{e;nor1}) \\

Thus,\\

\begin{align}
 & \ominus Hom(\Omega^{1}_{\CP^{1}}(a_{1} + \ldots + a_{k}),
\mathcal{O}_{\CP^{1}}) \oplus H^{0}(\CP^{1}, f^{*}(T \CP^{2})) \notag \\
& \oplus H^{1}(\CP^{1}, T \CP^{1}(-a_{1} - \ldots - a_{k})) \notag \\
\cong & [\ominus H^{0}( \CP^1, T \CP^1)
 \oplus ( T_{a_{1}}\CP^{1} \oplus \ldots \oplus T_{a_{k}} \CP^1) \notag \\
& \ominus H^{1}( \CP^1, T \CP^{1}( -a_{1} - \ldots - a_{k}))]
\bigoplus [ H^{0}( \CP^1, T \CP^1) \oplus H^{0}( \CP^1, N_{f}) ] \notag \\
& \bigoplus H^{1}(\CP^{1}, T
\CP^{1}(-a_{1} - \ldots - a_{k})) \notag \\
 \cong & H^{0}(\CP^1, N_{f}) \oplus T_{a_{1}} \CP^1 \oplus \ldots
\oplus T_{a_{k}} \CP^{1} \notag
\end{align}

\hfill q.e.d.

\vspace{2mm}

\begin{remark} \label{r;tricase} Let $\overline{M}_{k}$ be the
Delign-Mumford moduli space of genus zero curves with $k$-marked
points. Then,  $\overline{M}_{k}(\CP^2,0)$ is isomorphic to
$\overline{M}_{k} \times \CP^2$. Thus, we may consider its tangent
space at $[(f, \CP^1, a_{1}, \ldots, a_{k})]$ is $H^{1}(\CP^1, T
\CP^1(-a_{1}- \ldots - a_{k})) \oplus T_{f(\CP^1)} \CP^2$. This
tangent space formula is identical to the formula in Lemma
\ref{l;tan} in a $K$-theoretic point of view:

\vspace{2mm}

$H^{0}(\CP^1, N_{f}) \oplus T_{a_{1}} \CP^1
\oplus \ldots \oplus T_{a_{k}} \CP^{1}$ \\
 $\cong  H^{0}(\CP^1, f^{*}T\CP^{2})\ominus H^{0}(\CP^1, T \CP^1) \oplus T_{a_{1}} \CP^1
\oplus \ldots \oplus T_{a_{k}} \CP^{1}$ by (\ref{e;nor}) \\
$ \cong H^{0}(\CP^1, f^{*}T\CP^{2}) \oplus H^{1}(\CP^1, T
\CP^1(-a_{1}- \ldots - a_{k}))$ by (\ref{e;sh}) since
$k$ is greater than 2 and the degree of $T\CP^1$ is 2 \\
$\cong T_{f(\CP^1)} \CP^2 \oplus H^{1}(\CP^1, T \CP^1(-a_{1}- \ldots
- a_{k}))$ because deg$f$=0
\end{remark}

\par \emph{The evaluation map $ev$ on $\md{k}{X}$} is a morphism from
$\md{k}{X}$ to $X \times \ldots \times X$ which sends $[(f, C,
a_{1}, \ldots, a_{k})]$ to $(f(a_{1}), \ldots, f(a_{k}))$. The
following Proposition also establishes a well-known transversality
property of intersection cycles on the nodal stable maps locus in
$I_{d}([point], \ldots, [point])$ case.

\begin{proposition}\label{p;trans}
The evaluation map $ev$ on $\md{3d-1}{\CP^2}$ is a local isomorphism
on the nodal stable maps locus.
\end{proposition}
\vspace{2mm}

Proof. Let $\textbf{c}:=[(f, \CP^{1}, a_{1}, \ldots, a_{3d-1})]$
represent a nodal stable map. Then, a differential map $dev$ at
$\textbf{c}$  is\\

\noindent $\begin{array}{ccccc} T_{\textbf{c}} \md{3d-1}{\CP^2}
\cong & H^{0}(\CP^{1}, N_{f}) & \oplus T_{a_{1}} \CP^{1}&
\oplus \ldots & \oplus T_{a_{3d-1}} \CP^{1} \\
& (s,&  v_{1},&  \ldots, & v_{3d-1})
\end{array}$

\begin{flushright}
$\begin{array}{cccc}
 \stackrel{ dev
\mid_{\textbf{c}}}{\longrightarrow} & T_{f(a_{1})} \CP^{2} & \times
\ldots \times &  T_{f(a_{3d-1})}
\CP^{2}\\
 \mapsto & (s \mid_{a_{1}} + df
\mid_{a_{1}}(v_{1}),&  \ldots, & s \mid_{a_{3d-1}} + df
\mid_{a_{3d-1}}(v_{3d-1}))
\end{array}$
\end{flushright}

\vspace{2mm}

\noindent where $N_{f}$ is coker( $df: T \CP^{1} \rightarrow f^{*} T
\CP^{2}$). Suppose that $dev \mid_{\textbf{c}}(s, v_{1}, \ldots,
v_{3d-1}) = (0, \ldots, 0)$. Then, $s \mid_{a_{1}} = \ldots = s
\mid_{a_{3d-1}} = 0$ and $df \mid_{a_{1}}(v_{1}) = \ldots = df
\mid_{a_{3d-1}}(v_{3d-1}) = 0$ because the vectors $s \mid_{a_{i}}$
and $df \mid_{a_{i}}(v_{i})$ are linearly independent. Since the
degree of the line bundle $N_{f}$ is $3d-2$, the global section
which vanishes at $3d-1$ points represents a trivial element in
$H^{0}(\CP^{1}, N_{f}) $. Moreover, $v_{i}$, $i= 1, \ldots, 3d-1$,
are zero because $f$ is an immersion. This proves that $ev$ is
injective. It also implies that $ev$ is surjective because the
domain and the target of the linear map $d ev \mid_{\textbf{c}}$
have the same dimension. \hfill q.e.d.

\vspace{2mm}

\begin{lemma}\label{l;index} Assume $p_{1}, \ldots, p_{3d-1}$
represent trivial cycles in $\CP^2$ in general position. Let
$\textbf{c}:=[(f, \CP^{1}, a_{1}, \ldots, a_{3d-1})]$ be a point in
$ev_{1}^{-1}(p_{1}) \cap \ldots \cap ev_{3d-1}^{-1}(p_{3d-1})$.
Suppose that the degree of the evaluation map  $ev: \md{3d-1}{\CP^2}
\rightarrow \CP^2 \times \ldots \times \CP^2$ at $\textbf{c}$ is
$k$. Then, the intersection multiplicity at $\textbf{c}$ is $k$.
\end{lemma}

Proof. Since the transversality property is an open condition, the
Lemma follows from the definition of the local degree of the map and
the intersection multiplicity at $\textbf{c}$. \hfill q.e.d.\\

\par If the evaluation map is regular at $\textbf{c}$ in Lemma
\ref{l;index}, then Lemma \ref{l;index} implies that the cycles
$ev_{1}^{-1}(p_{1}), \ldots,  ev_{3d-1}^{-1}(p_{3d-1})$ meet
transversally at $\textbf{c}$.

\begin{proposition}\label{p;multone}
Assume $p_{1}, \ldots, p_{3d-1}$ represent trivial cycles in $\CP^2$
in general position. Let $\textbf{c}:=[(f, \CP^{1}, a_{1}, \ldots,
a_{3d-1})]$ be a point in $ev_{1}^{-1}(p_{1}) \cap \ldots \cap
ev_{3d-1}^{-1}(p_{3d-1})$ which represents either a nodal stable
map, a tacnode stable map or a triple node stable map. Then, the
intersection multiplicity at $\textbf{c}$ is one.
\end{proposition}

Proof. The case of a nodal stable map follows from Proposition
\ref{p;trans} and Lemma \ref{l;index}. A tacnode stable map is
gotten by the normalization of the image curve twice. Each
normalization is an immersion. Thus, a tacnode stable map is an
immersion. A triple node stable map is the normalization of the
image curve, which is an immersion. In either cases, $N_{f}$ is a
line bundle of degree $3d-2$. Thus, one can follow the proof of
Proposition \ref{p;trans} to show that $ev$ is an immersion at
$\textbf{c}$ if $\textbf{c}$ represents either a tacnode stable map
or a triple node stable map. The Proposition follows from Lemma
\ref{l;index}. \hfill q.e.d.

\begin{remark}\label{r;tac}
Despite of the results in Proposition \ref{p;multone}, the tacnode
stable maps locus, the triple node stable maps locus are not the
loci on which the intersection theory for an enumerative implication
is established. The reason is those loci are codimension one loci in
$\md{3d-1}{\CP^2}$.
\end{remark}

We calculate the kernel and the cokernel of an evaluation map $ev$
on $\md{3d-1}{\CP^2}$.

\begin{lemma}\label{l;kerofdev}
Let $\textbf{c}:= [(f,\CP^{1},a_{1}, \ldots,a_{3d-1})] \in
\md{3d-1}{\CP^2}$ represent a stable map such that $df \mid_{a_{i}}
\neq 0$, $i=1, \ldots, 3d-1$. Then, the kernel of the differential
of the $ev$ map at $\textbf{c}$ is isomorphic to $ H^{0}(\CP^1,
N(-a_{1}- \ldots - a_{3d-1}))$, where $N$ is coker$( df: T \CP^1
\rightarrow f^{*} T \CP^2)$.
\end{lemma}

\vspace{2mm}

Proof. Recall that $dev$ at $\textbf{c}$  is\\

\noindent $\begin{array}{ccccc} T_{\textbf{c}} \md{3d-1}{\CP^2}
\cong &  H^{0}(\CP^{1}, N) & \oplus T_{a_{1}} \CP^{1}&
\oplus \ldots & \oplus T_{a_{3d-1}} \CP^{1} \\
& (s,&  v_{1},&  \ldots, & v_{3d-1})
\end{array}$

\begin{flushright}
$\begin{array}{cccc}
 \stackrel{ dev
\mid_{\textbf{c}}}{\longrightarrow} & T_{f(a_{1})} \CP^{2} & \times
\ldots \times &  T_{f(a_{3d-1})}
\CP^{2}\\
 \mapsto & (s \mid_{a_{1}} + df
\mid_{a_{1}}(v_{1}),&  \ldots, & s \mid_{a_{3d-1}} + df
\mid_{a_{3d-1}}(v_{3d-1}))
\end{array}$
\end{flushright}

\vspace{2mm}

Since $f$ is a local immersion at $a_{i}$, $df \mid_{a_{i}}(v_{i}) =
0$ only if $v_{i} = 0$. $s \mid_{a_{i}} $, $i=1, \ldots, 3d-1$,
vanish only if $s \in H^{0}(\CP^1, N(-a_{1}-\ldots -a_{3d-1}))$. $df
\mid_{a_{i}}(v_{i})$ and $s \mid_{a_{i}}$ are independent vectors if
both are non-trivial. Thus, the kernel of $dev$ is $H^{0}(\CP^1,
N(-a_{1}-\ldots -a_{3d-1})) \oplus 0 \oplus \ldots \oplus 0 \cong
H^{0}(\CP^1, N(-a_{1}-\ldots -a_{3d-1}))$. \hfill q.e.d.

\vspace{2mm}

Let $N$ be coker$( df: T \CP^1 \rightarrow f^{*} T \CP^2)$. If $f$
has $k$ singularities, then $N$ has $k$ skyscraper sheaves which are
supported by critical points. The complement of $k$ skyscraper
sheaves in $N$ is a locally free sheaf. We will denote this locally
free sheaf by $NB_{k}$ and will call it \emph{the normal bundle of
$f$.}

\begin{lemma}\label{l;cokerofdev}Let $\textbf{c}:= [(f,\CP^{1},a_{1}, \ldots,a_{3d-1})] \in
\md{3d-1}{\CP^2}$ represent a stable map $f$ such that $df
\mid_{a_{i}} \neq 0$, $i=1, \ldots, 3d-1$. Suppose that $f$ has
exactly $k$ singular points $b_{1}, \ldots, b_{k}$ of degrees
$d_{1}, \ldots, d_{k}$ respectively. Then, the cokernel of the
differential of the $ev$ map at $\textbf{c}$ is isomorphic to $
H^{1}(\CP^1, N(-a_{1}- \ldots - a_{3d-1}))$ $\cong$ $H^{1}(\CP^1,
NB_{k}(-a_{1}- \ldots - a_{3d-1}))$.
\end{lemma}

\vspace{2mm}

Proof. We have the following short exact sequences of sheaves:

\begin{gather}
0 \rightarrow T \CP^1 \rightarrow f^{*} T \CP^2 \rightarrow N
\rightarrow 0 \label{e;nor} \\
0 \rightarrow \oplus_{i=1}^{k} \hspace{1mm} \tau_{i}^{d_{i}-1}
\rightarrow N \rightarrow NB_{k} \rightarrow 0 \label{e;sk}
\end{gather}

Lemma \ref{l;tan} and (\ref{e;sk}) shows that the tangent space at
$\textbf{c}$ is

\[ T_{\textbf{c}}\md{3d-1}{\CP^2} \cong \bigoplus_{i=1}^{k}
\tau_{i}^{d_{i}-1} \bigoplus H^{0}(\CP^1, NB_{k}) \bigoplus
\bigoplus_{i=1}^{3d-1} T_{a_{i}} \CP^1 \]

Lemma \ref{l;kerofdev} and (\ref{e;sk}) implies that the kernel of
$dev$ is \[ H^{0}(\CP^1, N(-a_{1}- \ldots - a_{3d-1})) \cong
\displaystyle \bigoplus_{i=1}^{k} \tau_{i}^{d_{i}-1}.\]

Thus, \\
\[ \mbox{coker} dev  \cong  \frac{ \bigoplus_{i=1}^{3d-1} T_{f(a_{i})}
\CP^2 }{dev( H^{0}( \CP^1, NB_{k})) \oplus dev(
\bigoplus_{i=1}^{3d-1} T_{a_{i}} \CP^1)}\]

Let $B$ and $C$ be subvector spaces of the vector space $A$ such
that $B \cap C = \{ 0 \}$. Then, one can easily check the elementary
isomorphism $ \displaystyle \frac{A}{ B \oplus C } \cong
\frac{C^{\perp} \oplus C }{B \oplus C} \cong \frac{C^{\perp}}{B} $,
where $C^{\perp}$ is an orthogonal complement of $C$ in $A$. Let
$N_{a_{i}}$ be the orthogonal complement of $dev(T_{a_{i}}\CP^1)$.
Since $dev( H^{0}( \CP^1, NB_{k})) \cap dev( \bigoplus_{i=1}^{3d-1}
T_{a_{i}} \CP^1) \cong \{ 0 \} $, we get

\begin{align}
 \mbox{coker} dev & \cong (\bigoplus_{i=1}^{3d-1} N_{a_{i}}) / dev(H^{0}( \CP^1,
 NB_{k}))\\
& \cong H^{0}(\CP^1, \bigoplus_{i=1}^{3d-1} \nu_{i})/ dev(H^{0}(
\CP^1, NB_{k})), \label{e;qu}
\end{align}

where $\nu_{i}$ is a skyscraper sheaf $N_{a_{i}}$ supported by
$a_{i}$, $i=1, \ldots, 3d-1$. Consider the short exact sequence of
sheaves $0 \rightarrow NB_{k}(-a_{1}- \ldots - a_{3d-1}) \rightarrow
NB_{k} \rightarrow  \bigoplus_{i=1}^{3d-1} \nu_{i} \rightarrow 0$.
Let's take a long exact sequence of sheaf cohomologies. Then, we
have

\begin{gather*}
0 \rightarrow H^{0}(\CP^1, NB_{k}(-a_{1}- \ldots - a_{3d-1}))
\rightarrow H^{0}(\CP^1, NB_{k}) \rightarrow \\
\rightarrow H^{0}(\CP^1, \bigoplus_{i=1}^{3d-1} \nu_{i} )
\rightarrow H^{1}( \CP^1, NB_{k}(-a_{1}- \ldots - a_{3d-1}))
\rightarrow \\ \rightarrow H^{1}(\CP^{1}, NB_{k})  \rightarrow
\ldots.
\end{gather*}

Note that

\begin{itemize}
\item $H^{0}(\CP^1, NB_{k}(-a_{1} - \ldots - a_{3d-1}))$ vanishes
because the degree of $NB_{k}$ is less than $3d-1$. \item The long
exact sequence induced by an exact sequence of sheaves: \[ 0
\rightarrow T \CP^{1} \rightarrow f^{*} T \CP^2 \rightarrow TN
\rightarrow 0 \] shows that $H^{1}(\CP^1, N)$ vanishes because
$\CP^2$ is a convex variety, i.e., $H^{1}(\CP^1, f^{*}T \CP^2)=0$.
Since $H^{1}(\CP^1, N)$ is isomorphic to $H^{1}(\CP^1, NB_{k})$,
$H^{1}(\CP^1, NB_{k})$ vanishes.
\end{itemize}

Thus, we get the desired isomorphism. \hfill q.e.d.\\

\begin{proposition}\label{p;cusp}
The cuspidal stable maps locus forms a degree 2 critical points set
of the evaluation map $ev$.
\end{proposition}
\vspace{2mm}

Proof. Let $\textbf{c} := [(f, \CP^1, a_{1}, \ldots, a_{3d-1})]$
represent a cuspidal stable map. Lemma \ref{l;kerofdev} implies the
kernel of $dev_{\textbf{c}} : = \tau$, where $\tau$ is a skyscraper
sheaf supported by the cuspidal singularity. cf. (\ref{e;sk}). Lemma
\ref{l;cokerofdev} implies the cokernel of $dev := H^{1}(\CP^1,
NB_{k}(-a_{1}- \ldots - a_{3d-1}))$. The Kodaira-Serre duality shows
that $H^{1}(\CP^1, NB(-a_{1}- \ldots - a_{3d-1}))$ is isomorphic to
$H^{0}(\CP^1, \omega_{\CP^1} \otimes [NB(-a_{1}- \ldots -
a_{3d-1})]^{*})$, where $\omega_{\CP^1}$ is a dualizing sheaf on
$\CP^1$ and $[NB(-a_{1}- \ldots - a_{3d-1})]^{*}$ is a dual vector
bundle of $[NB(-a_{1}- \ldots - a_{3d-1})]$. There is a canonical
residue morphism of degree 2 from the local slice of the direction
$\tau$ to the direction  $H^{1}(\CP^1, NB_{k}(-a_{1}- \ldots -
a_{3d-1}))$ induced by the $ev$ map. The morphism is:

\begin{align}
 \tau \simeq \tau^{*} &
\rightarrow & H^{0}(\CP^1, \omega_{\CP^1} \otimes [NB(-a_{1}- \ldots
-
a_{3d-1})]^{*})^{*} \notag \\
v dz & \mapsto & \frac{1}{2 \pi i}\int_{\gamma} \frac{v^{2}}{z} dz
\hspace{5cm}  \notag
\end{align}
because the local index of $f$ is 2. Thus, the Proposition follows.
 \hfill q.e.d. \\

\begin{remark} \label{r;dimc}
By the results in this section, we get the following equivalent
conditions:
\begin{itemize}
\item $dev$ has a critical point at $[(f, \CP^1,a_{1}, \ldots,
a_{3d-1})]$. \item $f$ is not an immersion. \item
$H^{0}(\CP^1,N(-a_{1} - \ldots -a_{3d-1}))$ doesn't vanish.
\item $H^{1}(\CP^1,N(-a_{1} - \ldots -a_{3d-1}))$ doesn't vanish,
\end{itemize}
\noindent where $N$ is Coker$(df: T \CP^1 \rightarrow f^{*} T
\CP^2)$.
\end{remark}

\vspace{2mm}

 Let $C$ be a pointed reducible curve which has two, pointed irreducible
components $(C_{1}, z_{1}, \ldots, z_{r})$, $(C_{2}, w_{1}, \ldots,
w_{s})$. Let $q_{1} \in C_{1}$, $q_{2} \in C_{2}$ be
(pre)gluing points.\\
 Then, a pointed stable map
$(f,C,z_{1}, \ldots, z_{r}, w_{1}, \ldots, w_{s})$ can be written as
$((f_{1}, (C_{1}, q_{1}), z_{1}, \ldots, z_{r}), (f_{2}, (C_{2},
q_{2}), w_{1}, \ldots, w_{s}))$.

\vspace{2mm}

\begin{lemma}\label{l;redt}
The tangent space splitting at \\
$[(f,C,z_{1}, \ldots, z_{r}, w_{1}, \ldots, w_{s})]:= [((f_{1},
(C_{1}, q_{1}), z_{1}, \ldots, z_{r}), (f_{2}, (C_{2}, q_{2}),
w_{1}, \ldots, w_{s}))]$ is

\vspace{2mm}

\begin{center}
$H^{0}(C_{1}, N_{1}) \oplus H^{0}(C_{2}, N_{2}) \oplus
T_{z_{1}}C_{1} \oplus \ldots \oplus T_{z_{r}}C_{1} \oplus
T_{w_{1}}C_{2} \oplus \ldots \oplus T_{w_{s}}C_{2} \oplus
(T_{q_{1}}C_{1} \otimes T_{q_{2}}C_{2})$ $\oplus T_{q_{1}}C_{1}
\oplus T_{q_{2}}C_{2} \ominus T_{f(q)} \CP^2$
\end{center}
\vspace{2mm}

\noindent where $N_{i}$ is a Coker($df_{i}: T C_{i} \rightarrow
f^{*} T \CP^2$), $i= 1,2$, and $T_{f(q)} \CP^2$ is a skyscraper
sheaf supported at $f(q)$, $q$ is a node in $C$, the degree of
$f_{i}$, $i=1,2$, is non-trivial.
\end{lemma}
\vspace{2mm}

Proof. We have to repeat the similar calculations we have done in
Lemma ~\ref{l;tan}. \\
From \hfill the \hfill long \hfill exact \hfill sequence \hfill
associated \hfill to \hfill the \hfill
hyperext \hfill group \\
$Ext^{1}(f^{*}\Omega^{1}_{\CP^2} \rightarrow \Omega^{1}_{C}(z_{1}+
\ldots + w_{s}), \mathcal{O}_{C})$:

\begin{multline*}
0 \hspace{1mm} \rightarrow \hspace{1mm} Hom(\Omega_{C}^{1}(z_{1} +
\ldots + z_{r} + w_{1} + \ldots + w_{s}), \mathcal{O}_{C})
\hspace{1mm} \rightarrow \hspace{1mm}
H^{0}(C, f^{*}T \CP^{2}) \rightarrow\\
\rightarrow Ext^{1}(f^{*}\Omega^{1}_{\CP^{2}} \rightarrow
\Omega^{1}_{C}(z_{1} + \ldots + z_{r} + w_{1} + \ldots + w_{s}),
\mathcal{O}_{C}) \rightarrow \\
\rightarrow Ext^{1}(\Omega^{1}_{C}(z_{1} + \ldots + z_{r} + w_{1} +
\ldots + w_{s}),\mathcal{O}_{C}) \rightarrow 0,
\end{multline*}

we get the following tangent space splitting at $(f,C,z_{1}, \ldots,
z_{r}, w_{1}, \ldots, w_{s})$:

\vspace{2mm}

\begin{gather}
\ominus Hom(\Omega_{C}^{1}(z_{1} + \ldots + z_{r} + w_{1} + \ldots
+ w_{s}), \mathcal{O}_{C}) \hspace{1mm} \oplus \label{e;all} \\
\oplus H^{0}(C, f^{*}T \CP^{2}) \hspace{1mm} \oplus \hspace{1mm}
Ext^{1}(\Omega^{1}_{C}(z_{1} + \ldots + z_{r} + w_{1} + \ldots +
w_{s}),\mathcal{O}_{C}). \notag
\end{gather}

A standard fact we will use in the following calculations is
$Hom(\Omega_{C}^{1}, \mathcal{O}_{C})$,
$\underline{Ext}^{0}(\Omega_{C}, \mathcal{O}_{C})$ are the sheaf of
derivations gotten by the pushforward of the sheaf of vector fields
on $\widetilde{C}:= C_{1} \cup C_{2}$ vanishing at the inverse
images $q_{1}, q_{2}$ of the node in $C$. Let $\pi: \widetilde{C}
 \rightarrow C$ be a normalization map.

\vspace{2mm}

\par We calculate the splitting of each term first.

\vspace{2mm}

\par For $\ominus Hom(\Omega_{C}^{1}(z_{1} + \ldots + z_{r} + w_{1} +
\ldots + w_{s}), \mathcal{O}_{C})$ term,\\
we use the short exact sequences of sheaves:
\begin{gather}
0 \rightarrow TC_{1}(-z_{1} - \ldots - z_{r}-q_{1}) \rightarrow
TC_{1} \rightarrow T_{z_{1}}C_{1} \oplus \ldots  \oplus
T_{z_{r}}C_{1} \oplus T_{q_{1}}C_{1} \rightarrow 0 \label{e;hom1}
\\
0 \rightarrow TC_{2}(-w_{1} - \ldots - w_{s}-q_{2}) \rightarrow
TC_{2} \rightarrow T_{w_{1}}C_{2} \oplus \ldots  \oplus
T_{w_{s}}C_{2} \oplus T_{q_{2}}C_{2} \rightarrow 0 \label{e;hom2}
\end{gather}

to get the following K-group equation  \vspace{2mm}

\begin{align}
 & Hom(\Omega_{C}^{1}(z_{1} + \ldots + z_{r} + w_{1} + \ldots +
w_{s}), \mathcal{O}_{C}) \notag \\
= & H^{0}(C, TC(-z_{1} - \ldots - w_{s})) \notag\\
= & H^{0}(C, \pi_{*}( T\widetilde{C}(-z_{1} - \ldots -w_{s}-q_{1}
-q_{2}))) \notag\\
= & H^{0}(\widetilde{C},  T \widetilde{C}(-z_{1} - \ldots
-w_{s}-q_{1} -q_{2})) \notag\\
= & H^{0}(C_{1}, TC_{1}(-z_{1} - \ldots - z_{r} - q_{1})) \oplus
H^{0}(C_{2}, TC_{2}(-w_{1} - \ldots - w_{s} - q_{2}))
\notag \\
= & H^{0}(C_{1}, TC_{1}) \ominus T_{z_{1}}C_{1}
\ominus \ldots \ominus T_{z_{r}}C_{1} \ominus T_{q_{1}}C_{1} \oplus
H^{1}(C_{1}, TC_{1}(-z_{1} - \ldots - z_{r} - \notag \\
& q_{1})) \oplus H^{0}(C_{2}, TC_{2})\ominus T_{w_{1}}C_{2} \ominus
\ldots \ominus T_{w_{s}}C_{2}\ominus T_{q_{2}} C_{2}\oplus
H^{1}(C_{2}, TC_{2}(-w_{1} \notag \\
& - \ldots - w_{s} - q_{2})) \notag
\end{align}

\noindent by (\ref{e;hom1}), (\ref{e;hom2}).

\vspace{2mm}

\par For $H^{0}(C, f^{*}T \CP^{2})$, we use the short exact sequence of
sheaves

\begin{gather*} 0 \rightarrow f^* T \CP^2
\rightarrow f_{1}^{*} T \CP^2 \oplus f_{2}^{*} T \CP^2 \rightarrow
T_{f(q)} \CP^2 \rightarrow 0,
\end{gather*}

\noindent to get a K-group equation

\begin{gather*}
H^{0}(C, f^{*}T \CP^{2}) = H^{0}(C_{1}, f_{1}^{*} T \CP^2) \oplus
H^{0}(C_{2}, f_{2}^{*} T \CP^2) \ominus T_{f(q)} \CP^2,
\end{gather*}

\noindent because $H^{1}(C, f^{*}T \CP^{2})$ vanishes by Lemma 10 in
~\cite{fap}.

\vspace{2mm}

\par For $Ext^{1}(\Omega^{1}_{C}(z_{1} + \ldots + z_{r} + w_{1} +
\ldots + w_{s}),\mathcal{O}_{C})$, we use an exact sequence from the
local to global spectral sequence in Lemma \ref{l;tan} to get

\vspace{2mm}

$Ext^{1}(\Omega^{1}_{C}(z_{1} + \ldots + z_{r} + w_{1} + \ldots +
w_{s}),\mathcal{O}_{C})$ \vspace{2mm}

\noindent $= H^{1}(C, \underline{Ext}^{0}(\Omega_{C}(z_{1}+ \ldots +
w_{s}), \mathcal{O}_{C})) \oplus H^{0}(C,
\underline{Ext}^{1}(\Omega_{C}(z_{1} + \ldots + w_{s}),
\mathcal{O}_{C}))$ \vspace{2mm}

\noindent $= H^{1}(C, \pi_{*}( T \widetilde{C}(-z_{1} - \ldots
-w_{s}-q_{1} -q_{2})))
 \oplus H^{0}(C, \underline{Ext}^{1}(\Omega_{C}(z_{1}
+ \ldots + w_{s}), \mathcal{O}_{C})) $ \vspace{2mm}

\noindent $= H^{1}(C_{1}, TC_{1}(-z_{1} - \ldots - z_{r}-q_{1}))
\oplus H^{1}(C_{2}, TC_{2}(-w_{1} - \ldots - w_{s}-q_{2})) \oplus
T_{q_{1}} C_{1} \otimes T_{q_{2}} C_{2}$.

\vspace{2mm}

\par The Lemma follows by putting all terms to (\ref{e;all}).
\hfill q.e.d.

\vspace{2mm}

\begin{remark}
1. Let $f_{2}$ be a degree 0 map. Similar calculations we have done
in Remark \ref{r;tricase} by using Lemma \ref{l;redt} shows that the
tangent space splitting at
\begin{multline*}
[(f,C,z_{1}, \ldots, z_{r}, w_{1}, \ldots, w_{s})]\\
:= [((f_{1},
(C_{1}, q_{1}), z_{1}, \ldots, z_{r}), (f_{2}, (C_{2}, q_{2}),
w_{1}, \ldots, w_{s}))],
\end{multline*}
is:

\vspace{2mm}

\begin{center}
$H^{0}(C_{1}, N_{1}) \oplus T_{z_{1}}C_{1} \oplus \ldots \oplus
T_{z_{r}}C_{1} \oplus T_{q_{1}}C_{1} \oplus$ $H^{1}(C_{2}, T
C_{2}(-w_{1} - \ldots - w_{s} -q_{2}))$ $\oplus (T_{q_{1}}C_{1}
\otimes T_{q_{2}}C_{2})$
\end{center}

\vspace{2mm}

2. One can extend the result in Lemma \ref{l;redt} to the general
case. Let $\mathbf{c}:=[(f,C, a_{1}, \ldots, a_{k})]$ be a point in
$\md{k}{\CP^2}$. Let $\pi : \widetilde{C} := \CP^{1}_{1} \cup \ldots
\cup \CP^{1}_{l} \rightarrow C$ be a normalization map, where
$\CP^{1}_{i}$ is biholomorphic to $\CP^{1}$. Let $g_{1}, \ldots,
g_{r}$ be singular points on $C$, $r:= l-1$. Let's denote elements
in $\pi^{-1}(g_{i})$ by $g_{i}^{1}, g_{i}^{2}$. Let $N_{i}$ be
coker$(df_{i}:T\CP^{1}_{i} \rightarrow T \CP^2)$, where $f_{i}:= f
\mid_{\CP^{1}_{i}}$. Then, the tangent space
$T_{\mathbf{c}}\md{k}{\CP^2}$ at $\mathbf{c}$ is
\begin{gather}
\bigoplus_{i=1}^{l} H^{0}(\CP^{1}_{i}, N_{i}) \oplus \bigoplus_{i=1,
\ldots,k} T_{a_{i}}\CP^{1}_{q(a_{i})} \oplus
(\bigoplus_{i=1,\ldots,r} T_{g_{i}^{1}}\CP^{1}_{q(g_{i}^{1})}\otimes
T_{g_{i}^{2}} \CP^{1}_{q(g_{i}^{2})}) \oplus   \notag \\
 \oplus \bigoplus_{i=1, \ldots,r}^{j=1,2}
T_{g_{i}^{j}}\CP^{1}_{q(g_{i}^{j})} \ominus(\bigoplus_{i=1}^{r}
T_{f(g_{i})}\CP^2). \notag
\end{gather}
See \cite{kwon} for details of calculations.
\end{remark}

\vspace{2mm}

The equivalent conditions in Remark \ref{r;dimc} are no longer true
if we consider reducible stable maps. The rank of the cokernel of an
evaluation map is determined by the number of marked points on each
irreducible domain curve and the mapping properties of the stable
map $f$ on each irreducible component.

\begin{proposition}\label{p;red}
Let
\begin{multline*}
\textbf{c} :=[(f,C,z_{1},\ldots, z_{r}, w_{1}, \ldots, w_{s})] \\
:= [((f_{1},
(C_{1}, q_{1}), z_{1}, \ldots, z_{r}), (f_{2}, (C_{2}, q_{2}),
w_{1}, \ldots, w_{s}))]
\end{multline*}
represent a reducible stable map, where $C_{i}$, $i=1,2$, is
isomorphic to $\CP^1$ and $f_{i}$ is an immersion of degree $d_{i}$,
$i = 1,2$ if $f$ is not trivial on the component. Then, \\
$(i)$ If any of $f_{i}$ is a degree 0 map, then the cokernel of
$dev$ at $\textbf{c}$ has a rank bigger than one. \\
$(ii)$ If $r$ or $s$ is strictly bigger than $3d_{1} -2$ or $3d_{2}
-2$ respectively, then the cokernel of $dev$ at $\textbf{c}$ has a
rank bigger
than two. \\
$(iii)$ If $r$ or $s$ is $3d_{1}-2$ or $3d_{2} -2$ respectively,
then the cokernel of $dev$ at $\textbf{c}$ has a rank one.\\
$(iv)$ If $r$ or $s$ is $3d_{1} -1$ or $3d_{2}-1$ respectively, then
the evaluation map $ev$ at $\textbf{c}$ is regular.
\end{proposition}
\vspace{2mm}

\begin{remark}
Bezout's theorem implies that any deformed image curves determined
by vectors in $H^{0}(C_{1}, N_{1})$ and $H^{0}(C_{2}, N_{2})$ always
meet if $f$ is non-trivial, when the target space dimension is two.
Thus, in this particular dimension, we can calculate the rank of the
cokernel of the $ev$ map by considering the vectors other than the
vectors in $V:= T_{q_{1}}C_{1} \oplus T_{q_{2}}C_{2} \ominus
T_{f(q)} \CP^2$.
\end{remark}

Sketch of the Proof of Proposition \ref{p;red}. \hspace{3mm} Suppose
that $f_{1}$ is a trivial map. The stability condition implies that
$C_{1}$ contains at least 2 marked points. The differential of the
$i$-th evaluation map $ev_{i}$ is zero on $C_{1}$. The maximum
dimensional contribution to the rank of the $ev$ map from $C_{2}$ is
at most $(3d-3) + (3d-3)$. (i) follows from this.\\

For (ii), (iii), (iv), we need to calculate the dimensional
contributions from $T_{q_{1}}C_{1} \otimes T_{q_{2}}C_{2}$  and

\begin{multline*} H^{0}(C_{1}, N_{1})
 \oplus
H^{0}(C_{2}, N_{2})  \oplus  T_{z_{1}}C_{1}   \oplus \ldots \oplus
T_{z_{r}}C_{1} \oplus
T_{w_{1}}C_{2}  \oplus \ldots \oplus  T_{w_{s}}C_{2} \\
\stackrel{dev}{\rightarrow}   T_{f(z_{1})} \CP^2 \hspace{3mm} \times
\ldots \times T_{f(z_{r})} \CP^2
\times  T_{f(w_{1})} \CP^2  \times \ldots \times  T_{f(w_{s})} \CP^2 \\
(s ,  t ,  v_{1},  \ldots,  v_{r}  , v'_{1} ,  \ldots,  v'_{s})
\mapsto  (s \mid_{z_{1}} + df \mid_{z_{1}}(v_{1}), \ldots,
 t \mid_{w_{s}} + df \mid_{w_{s}}(v'_{s}))
\end{multline*}

\vspace{2mm}

We use the following facts:
\begin{itemize}
\item $T_{q_{1}}C_{1} \otimes T_{q_{2}}C_{2}$ contributes to the rank
of $dev$ by one.
\item The vectors $s \mid_{z_{i}}$ and $df \mid_{z_{i}}(v_{i})$, $i=1,
\ldots, r$, $t \mid_{w_{j}}$ and $df \mid_{w_{j}}(v'_{j})$, $j=1,
\ldots, s$ are linearly independent.
\end{itemize}

The reason for the first item is $T_{q_{1}}C_{1} \otimes
T_{q_{2}}C_{2}$ generates the first order deformation of smoothing
node. During the smoothing node deformation, the stable map also
changes because the stable maps are the same if they agree to
infinite order at any point. Since the deformation space generating
smoothing node deformation is smooth, the image of any non-trivial
vector in $T_{q_{1}}C_{1} \otimes T_{q_{2}}C_{2}$ by $dev$ is
non-trivial.
\par The result follows from the straightforward dimension counts.
\hfill q.e.d.

\vspace{2mm}

Let $V_{1}, \ldots, V_{k}$ be the intersection cycles in the variety
$V$. Let $\textbf{p}$ be a point in $V_{1} \cap \ldots \cap V_{k}$.
We will say \emph{`the transversality property is not established at
$\textbf{p}$'} if the intersection multiplicity at $\textbf{p}$ is
not finite. One can easily check that the loci stated in Proposition
\ref{p;red} (i) - (iii) are examples of that. Note that the image of
the loci in Proposition \ref{p;red} (iii) is the subloci of
codimension 2 in $\CP^2 \times \ldots \times \CP^2$ in general
points, which is different from the rank of $dev$ along these loci.

\begin{theorem}\label{p;gent}
Let $\textbf{p}$ be a point in $ev_{1}^{-1}(p_{1}) \cap \ldots \cap
ev_{3d-1}^{-1}(p_{3d-1})$ which represents a cuspidal stable map,
where $p_{i}$ is a Chow 0-cycle's representatives,
 $i= 1, \ldots, 3d-1$.
Then, the intersection multiplicity at $\textbf{p}$ is 2. That is, a
transversality uniformly fails along the cuspidal stable maps locus.
A cuspidal stable maps locus is a unique equi-singular locus which
has a codimension $\leq 1$ on which a transversality always fails.
\end{theorem}
\vspace{2mm}

Proof. Proposition \ref{p;multone} implies that if $\textbf{p}$
represents a nodal stable map, a tacnode stable map, a triple node
stable map, then the intersection multiplicity at $\textbf{p}$ is
one. Proposition \ref{p;cusp} and Lemma \ref{l;index} implies that
if $\textbf{p}$ represents a cuspidal stable map, then the
intersection multiplicity at $\textbf{p}$ is 2. The classification
of the codimension $\leq 1$ equi-singular loci in Corollary
\ref{c;dense} and Theorem \ref{t;stdeg} implies the last statement
of the Theorem. \hfill q.e.d.

\vspace{2mm}

\subsubsection{Real Version of the Gromov-Witten invariants:\\
 $I_{d}([point], \ldots, [point])$ case} \label{ss;point}

\par If the target space $X$ is equipped with a real structure
$\tau$, then it induces a real structure $\overline{M}_{n}(X,
\beta)$ defined by $[(f, C, a_{1}, \ldots, a_{n})] \mapsto
[(\overline{f}, \overline{C}, \overline{a}_{1}, \ldots,
\overline{a}_{n})]$, where $\overline{f}(x) = \tau \circ f(
\overline{x})$, $\overline{C} = (C_{1}, \overline{q}_{1}) \cup
\ldots \cup (C_{m}, \overline{q}_{m})$ if $C = (C_{1}, q_{1}) \cup
\ldots \cup (C_{m}, q_{m})$ and $C_{i}$ is isomorphic to $\CP^1$.
The evaluation map $ev$ is a real map, that is, it commutes with the
real structure of the moduli space and of the target space. See
\cite{kwon}. Let's denote the real part of the moduli space by
$\overline{M}_{n}(X, \beta)^{re}$. Then, the $ev$ map sends
$\overline{M}_{n}(X, \beta)^{re}$ to the real part of the target space.\\

\par Theorem \ref{p;gent} doesn't have the enumerative
implications in the complex Gromov-Witten theory. However, Theorem
\ref{p;gent} is very important when we define the real version of
local invariants. The cuspidal stable maps locus comparts the moduli
space $\md{3d-1}{\CP^2}^{re}$. The image of the cuspidal stable maps
locus by the $ev$ map comparts the $3d-1$ fold product of $\RP^2$.
The notions of chambers and walls in the moduli space
$\md{3d-1}{\CP^2}^{re}$ and $\RP^{2} \times \ldots \times \RP^{2}$
arise.

\begin{Def}(Chambers and Walls in the moduli space
$\md{3d-1}{\CP^2}$)\\
\emph{Walls} in $\md{3d-1}{\CP^2}^{re}$ are the codimension one loci
on which transversality uniformly fails. \emph{Chambers} in
$\md{3d-1}{\CP^2}^{re}$ are the connected components in the
complement of walls.
\end{Def}

\begin{remark}
Theorem \ref{p;gent} implies that $\textbf{c}:=[(f, C, a_{1},
\ldots, a_{3d-1})]$ belongs to the wall if and only if $\textbf{c}$
represents a cuspidal stable map.
\end{remark}

We will call the $3d-1$ fold product $\RP^{2} \times \ldots \times
\RP^2$ of $\RP^2$ as a \emph{real Chow 0-cycles parameter space}.

\begin{Def}
\emph{Walls} in the real Chow 0-cycles parameter space are the
codimension one regions of the image of the cuspidal stable maps
locus in $\md{3d-1}{\CP^2}^{re}$ by the evaluation map $ev$.
\emph{Chambers} are the connected components of the complement of
walls.
\end{Def}

\begin{remark}
If $\textbf{p}:= (p_{1}, \ldots, p_{3d-1}) \in \RP^{2} \times \ldots
\times \RP^2 $ is a general point in a chamber, then each point in
$\textbf{c} \in ev^{-1}(\textbf{p}):= ev_{1}^{-1}(p_{1}) \cap \ldots
\cap ev_{3d-1}^{-1}(p_{3d-1})$ represents a nodal stable map.
$\textbf{p}$ is in the wall if one of the points in
$ev^{-1}(\textbf{p})$ represents a cuspidal stable map and all other
points represent nodal stable maps.
\end{remark}

Obviously, the number of points in $ev^{-1}(\textbf{p})$ doesn't
change for general points in a chamber. Therefore, the local version
of the real Gromov-Witten invariants in the following Definition is
well-defined.

\begin{Def} Let $\mathcal{C}$ be a chamber in the Chow 0-cycles
parameter space. \emph{The real Gromov-Witten invariant} of the
chamber $\mathcal{C}$ is the number of points in $ev_{1}^{-1}(p_{1})
\cap \ldots \cap ev_{3d-1}^{-1}(p_{3d-1}) \cap
\md{3d-1}{\CP^2}^{re}$.
\end{Def}

\par The following Corollary show the differences
of the real Gromov-Witten invariants in adjacent chambers are
exactly two. And walls in $\md{3d-1}{\CP^2}$, $\RP^2 \times \ldots
\times \RP^2$ are the place we gain or loose two real solutions.

\begin{corollary}\label{c;realtr}
Let $\textbf{c}:=[(f, \CP^1, a_{1}, \ldots, a_{3d-1})]$ be a point
in the wall in $\md{3d-1}{\CP^2}^{re}$. Let $v$ be a real vector in
$coker(d ev \mid_{\textbf{c}})^{re}$ and $p: [-1,1] \rightarrow
\RP^2 \times \ldots \times \RP^2$ be a path satisfying the
followings:

\vspace{2mm}

\begin{itemize}
\item The path $p$ is tangential to the vector $v$ at $p(0)$.
\item $p([-1,0))$ is sitting in one chamber and $p((0,1])$ is
sitting in an adjacent chamber.
\end{itemize}
\vspace{2mm}

Let $\mathcal{H}$ be a neighborhood of $\textbf{c} \in p^{-1}(t)
\cap \md{3d-1}{\CP^2}^{re} \cap \mathcal{H}$, where $\textbf{c}$
represents a cuspidal stable map. Then, $p^{-1}(t) \cap
\md{3d-1}{\CP^2}^{re} \cap \mathcal{H}$ consists of two elements for
any general point $t \in [-1,0)$ and $p^{-1}(t) \cap
\md{3d-1}{\CP^2}^{re} \cap \mathcal{H}$ is an empty set for any
general point $t \in (0,1]$, or vice versa.
\end{corollary}
\vspace{2mm}

Proof. Proposition \ref{p;cusp} shows that the local model of the
$ev$ map along the local slice perpendicular to the cuspidal stable
maps locus is $z^{2} -t$ or $z^{2} + t$. The result follows by
restricting ourselves to the real part of the moduli space
$\md{3d-1}{\CP^2}^{re}$. \hfill q.e.d.\\

\par When we cross the wall in $\RP^2 \times \ldots \times
\RP^2$, the inverse images by the $ev$ map vary within the chambers
in $\md{3d-1}{\CP^2}^{re}$ except one pair of points which vary in
adjacent chambers in $\md{3d-1}{\CP^2}^{re}$ and meet at the wall in
$\md{3d-1}{\CP^2}^{re}$ and disappear.

\par In real algebraic point of view, we have two notions of node singularities.

\vspace{2mm}

\begin{Def}
\emph{A non-isolated node} is a node whose local equation is $z^{2}
- w^{2} = 0$ or $zw=0$ by real coordinate changes. \emph{An isolated
node} is a node whose local equation is $z^{2} + w^{2}= 0$ by real
coordinate changes.
\end{Def}

\vspace{2mm}

The following example shows that refined sense's real singularities
induce different global topological invariants in the real part of
image curves. However, they don't have a topological implication in
complex curves.

\begin{example}
Let $T^{isol}$ be a rational real nodal curve of degree 3 which has
an isolated node and $T^{non-isol}$ be a rational real nodal curve
of degree 3 which has a non-isolated node. Then, they are isomorphic
in complex sense.
\par  Let $T$ be a real torus, i.e., a smooth torus
in $\CP^2$ represented by a degree 3 real polynomial. Then, $T$ has
two generators $\alpha$, $\beta$ in the fundamental group of $T$.
The self-automorphism sending $\alpha$ to $\beta$ and $\beta$ to
$\alpha$ induces an isomorphism between $T^{isol}$ and
$T^{non-isol}$ because if $T^{isol}$ is gotten by trivializing a
generator $\alpha$, then $T^{non-isol}$ is from trivializing a
generator $\beta$. The Euler characteristics $\chi(T^{isol})$,
$\chi(T^{non-isol})$ are both one. However, the Euler characteristic
of the real part of $T^{isol}$ is 1 and that of the real part of
$T^{non-isol}$ is $-1$.
\end{example}

\begin{figure}
\includegraphics{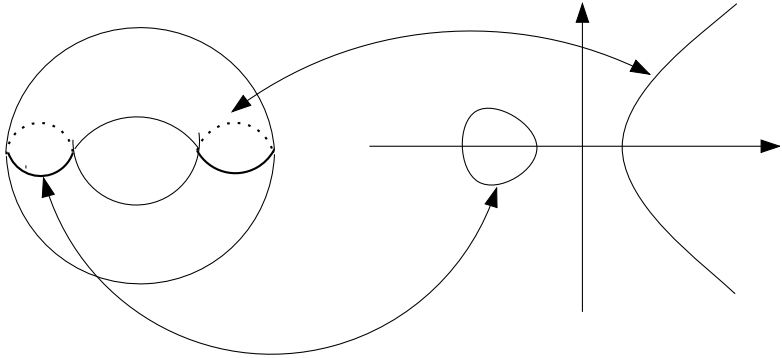}
\caption{Smooth torus} \label{p;smoothtorus}
\end{figure}

\vspace{2mm}

\begin{figure}
\includegraphics{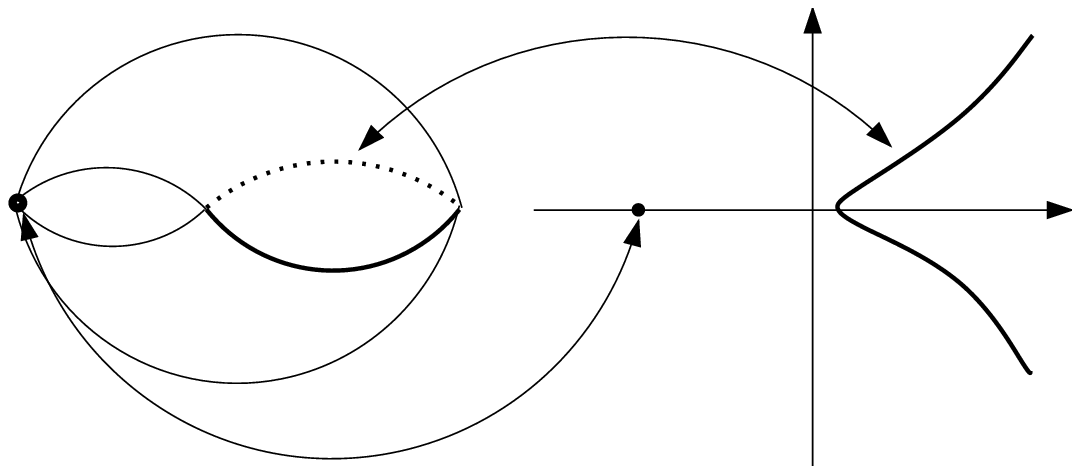}
\caption{Singular torus having an isolated node}
\label{p;isolatedtorus}
\end{figure}

\vspace{2mm}

\begin{figure}
\includegraphics{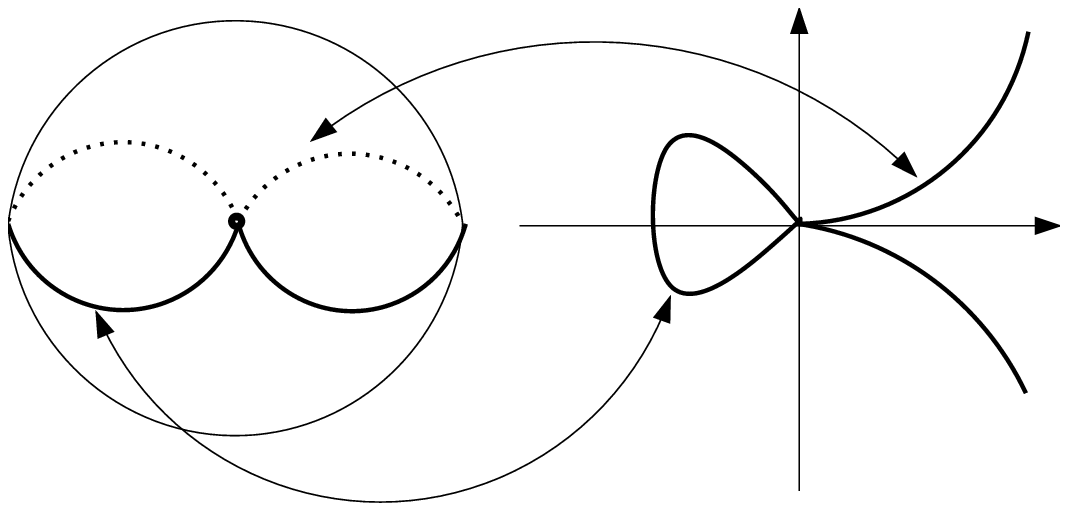}
\caption{Singular torus having a non-isolated node}
\label{p;nonisolatedtorus}
\end{figure}

\begin{figure}
\includegraphics{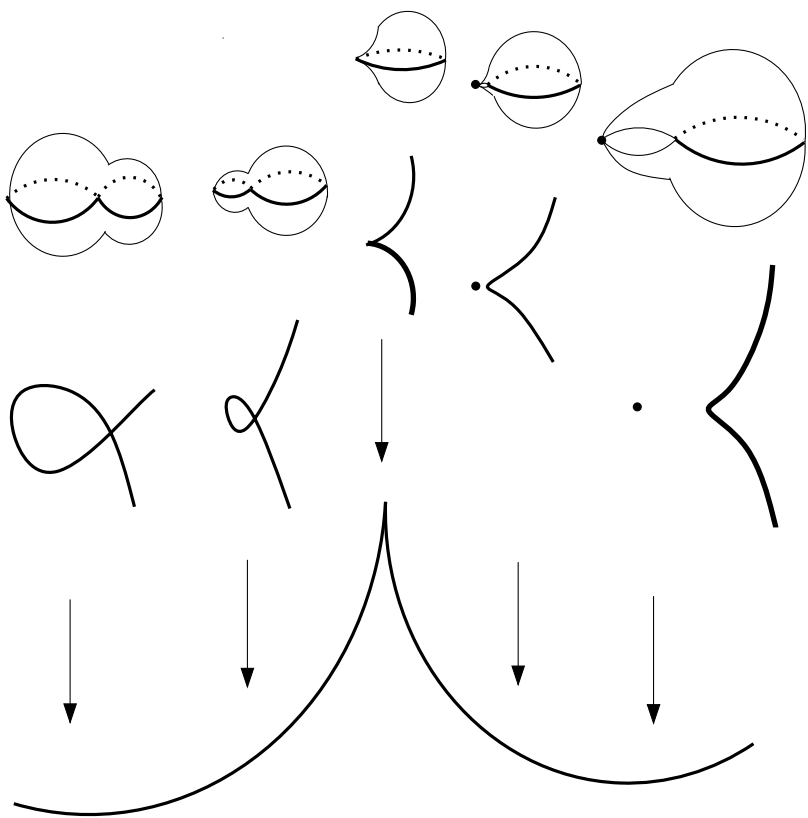}
\caption{Deformation of elliptic curves around a cuspidal curve}
\label{p;deformation}
\end{figure}

\begin{remark}\label{r;isono}
A cuspidal singularity is from the trivialization of two consecutive
generators in a smooth torus of a genus $g$. If trivializing one
generator created an isolated node, then trivializing the adjacent
generator created a non-isolated node. The topological transition in
the real part of the image curves happens only after we pass through
the cuspidal stable maps locus, that is, wall in the real part of
the moduli space $\md{3d-1}{\CP^2}$.  The difference in the number
of isolated nodes is exactly one. The number of isolated nodes and
the number of non-isolated nodes are topological invariants within a
chamber.
\end{remark}

There doesn't have any canonical way to give an orientation on the
real part of the moduli space $\md{3d-1}{\CP^2}$ because the real
part $\RP^2$ of the target space is non-orientable. Therefore, we
cannot count curves on the real part of the moduli space
$\md{3d-1}{\CP^2}$ intersection theoretically. Alternative approach
was taken by J-Y Welschinger in ~\cite{wel}. The Welschinger's
invariant provides the global minimum bound. The results we have
gotten so far show the invariance in an algebraic category when the
target space is $\RP^2$. However, in general target space cases, it
isn't necessary that every non-trasversality result produces minimum
bound results even for the $I_{\beta}([point], \ldots, [point])$
case.

\begin{proposition}\label{p;glo}
Let $ (p_{1}, \ldots, p_{3d-1})$ be any general element in $\RP^{2}
\times \ldots \times \RP^2$. Then, the real Gromov-Witten invariants
on any chamber is greater than or equal to

\vspace{2mm}

\begin{center}
$\mid$ $\displaystyle \sum_{k} (-1)^{k} \mbox{the number of real
stable maps in} \hspace{1mm} ev_{1}^{-1}(p_{1}) \cap \ldots \cap
ev_{3d-1}^{-1}(p_{3d-1}) \cap \md{3d-1}{\CP^{2}}^{re}$
 having $k$ isolated nodes $\mid$.
\end{center}
\end{proposition}

\vspace{2mm}

Proof. The result follows immediately from Corollary \ref{c;realtr}
and Remark \ref{r;isono} \hfill q.e.d.

\subsubsection{Transversality properties for the Gromov-Witten invariant
and its real enumerative implications:\\
\indent $I_{d}$([point], \ldots, [point],[line], \ldots, [line])
case}\label{ss;def}

\vspace{2mm}

The set of divisors linearly equivalent to $d$[line], $d >0$, is in
one-to-one correspondence with $(H^{0}(\CP^{2}, \mathcal{O}(d))
\setminus \{0 \})/ \C^{*}$, which is isomorphic to $\CP^{N(d)}$,
where $N(d) = \frac{d \cdot (d+3)}{2}$. We will denote it by
$\CP(d)$. The divisor class we will consider in this section will be
the combination of [point], [line], i.e., $d=$ 0 or 1. Let's denote
the product of the Chow cycles parameter space $\CP(d_{1}) \times
\ldots \times \CP(d_{k})$ by $\CP(d_{1}, \ldots, d_{k})$, where
$d_{i}=$ 0 or 1.

\begin{remark}\label{r;min}
Let $l$ be the number of non-zero $d_{i}$ and $m$ be the number of
trivial $d_{i}$ in $\CP(d_{1}, \ldots, d_{k})$.  The Gromov-Witten
invariant is enumerative, only when $2 \cdot m + l$ equals to the
dimension of $\md{k}{\CP^2}$. One can easily check that $3d-1$ is
the minimum number of marked points which produce the non-trivial
Gromov-Witten invariant. If $k > 3d-1$, then there are exactly
$3d-1$ number of $d_{i} =0$ and $k-(3d-1)$ number of $d_{j}=1$.
\end{remark}

\par Various transversality results in section \ref{ss;point}
are extendible as follows. However, the following Theorem doesn't
characterize the intersection theoretic properties as in the
$I_{d}$([point], \ldots, [point]) case.

\begin{theorem}\label{t;cusp}
Let the Gromov-Witten invariant $I_{d}(d_{1}, \ldots, d_{k})$,
$d_{i}=$ 0 or 1, be enumerative. Let $\textbf{p}$ be a point in
$ev_{1}^{-1}(\Lambda_{1}) \cap \ldots \cap
ev_{k}^{-1}(\Lambda_{k})$, where $\Lambda_{i}$ is a Chow 0- or 1-
cycle's representative,
$i= 1, \ldots, k$. \\
(i) Suppose that $\textbf{p}$ represents a cuspidal stable map and
the stable map represented by a point $\textbf{p}$ meets all Chow
1-cycle's representatives in $\{ \Lambda_{1}, \ldots, \Lambda_{k}
\}$ transversally. Then, the intersection multiplicity at
$\textbf{p}$ is 2. \\
(ii) Transversality uniformly fails along the
cuspidal stable maps
locus. \\
(iii) Suppose that $\textbf{p}$ represents a stable map which is
either a nodal stable map or a triple node stable map or a tac node
stable map. Assume that the stable map represented by a point
$\textbf{p}$ meets all Chow 1-cycle's representatives in $\{
\Lambda_{1}, \ldots, \Lambda_{k} \}$ transversally. Then, the
intersection multiplicity at $\textbf{p}$ is 1.
\end{theorem}

\vspace{2mm}

Proof. (i) By Remark ~\ref{r;min}, we may rearrange the cycles so
that we have $\CP(d_{1}, \ldots, d_{k}) = \CP(0, \ldots, 0, d_{3d},
\ldots, d_{k})$, $d_{i} =1 $ for $3d \leq i \leq k$.
\par Let $\textbf{p}:=[(g, \CP^1, a_{1}, \ldots, a_{k})]$ be a cuspidal stable map which is in \\
$ev_{1}^{-1}(\Lambda_{1}) \cap \ldots \cap ev_{k}^{-1}(\Lambda_{k})$
and satisfies the assumptions. Then,

\begin{itemize}
\item $g(\CP^1)$ passes through $\Lambda_{i}$ and $g(a_{i}) =
\Lambda_{i}$, $1 \leq i \leq 3d-1$.
\item $g(\CP^1)$ meets
$\Lambda_{j}$ transversally and $g(a_{j}) \in \Lambda_{j}$, $3d \leq
j \leq k$.
\end{itemize}
\vspace{2mm}

Clearly, the curve $g(\CP^1)$ meets the perturbed Chow 1-cycles
$\Lambda'_{j}$ of $\Lambda_{j}$, $3d \leq j \leq k$, transversally.
And $\textbf{p}$ varies in a unique way according to the deformation
of marked points whose image meet $\Lambda'_{j}$ if $\Lambda_{i}$,
$i=1, \ldots, 3d-1$, is fixed. Proposition \ref{p;cusp} and Lemma
\ref{l;index} show that if we perturb the points $\Lambda_{i}$, $1
\leq i \leq 3d-1$, then we get two nodal, pointed stable maps which
pass through the perturbed points. These pointed stable maps still
meet the Chow 1-cycles $\Lambda_{j}$, $3d \leq j \leq k$
transversally because the transversality property is an open
condition. Thus, the intersection multiplicity at $\textbf{p}$ is
two. (ii) is an immediate consequence of (i). \par (iii) Lemma
\ref{l;index} and Proposition \ref{p;multone} show that the
intersection multiplicity at $\textbf{p}$ is one because of the
transversality assumption in $(iii)$.
\hfill q.e.d.\\

\par In $I_{d}$([point], \ldots, [point]) case, the moduli space
$\md{3d-1}{\CP^2}$ and the Chow 0-cycles parameter space were
related by the $ev$ map. Wall and chamber notions in
$\md{3d-1}{\CP^2}^{re}$ and $\RP^2 \times \ldots \times \RP^2$ can
be characterized by topological and geometric significance. That is
no more the case if we consider the $I_{d}$([point], \ldots,
[point],[line], \ldots, [line]) case. We don't have a map which
relates a moduli space with the Chow cycles parameter space.
Regardless of Theorem \ref{t;cusp}, the characterization of the
intersection theoretic properties cannot be extended to the
$I_{d}$([point], \ldots, [point],[line], \ldots, [line]) case.
Nevertheless, there are notions of walls and chambers in the Chow
cycles parameter space.
\par Recall the following well-known real transversality principle:

\vspace{2mm}

\par Small real perturbations of a transverse intersection preserve
transversality as well as the number of real and complex points in
the intersection.

\vspace{2mm}

\par Based on the real transversality principle, we can extend the
notion of the real version of the Gromov-Witten invariants on the
cycles parameter space as follows:

\begin{Def}
\emph{Wall} $\mathcal{W}$ is a codimension one locus in the cycles
parameter space $\CP(d_{1}, \ldots, d_{k})^{re}$, $d_{i}=$ 0 for $i
=1, \ldots, 3d-1$ and $d_{j} = 1$ for $j > 3d-1$, such that all
elements except one in $ev_{1}^{-1}(\Lambda_{1}) \cap \ldots \cap
ev_{k}^{-1}(\Lambda_{k}) \cap \md{k}{\CP^2}^{re}$ have an
intersection multiplicity one, where $\Lambda_{1}, \ldots,
\Lambda_{k}$ are real ordered cycles represented by an element
$(\lambda_{1}, \ldots, \lambda_{k})$ in $\mathcal{W}$.
\par \emph{A Chamber} \hfill in \hfill $\CP(d_{1}, \ldots, d_{k})^{re}$ \hfill is \hfill a
\hfill connected \hfill component \hfill in \\
$\CP(d_{1}, \ldots, d_{k})^{re}\setminus \mathcal{W}$. \end{Def}

For a general element $(\lambda_{1}, \ldots, \lambda_{k})$ in the
chamber, the intersection multiplicity of any point in
$ev_{1}^{-1}(\Lambda_{1}) \cap \ldots \cap ev_{k}^{-1}(\Lambda_{k})
\cap \md{k}{\CP^2}^{re}$ is one. The following real version of the
Gromov-Witten invariant on the chamber is well-defined by the real
transversality principle:

\begin{Def}
Let \hfill $(\lambda_{1}, \ldots, \lambda_{k})$ \hfill be \hfill a
\hfill general \hfill element \hfill in \hfill a \hfill given \\
chamber $\mathcal{C} \subseteq \CP(d_{1}, \ldots, d_{k})^{re}$ which
represents the ordered cycles' representatives $\Lambda_{1}, \ldots,
\Lambda_{k}$, $k \geq 3d-1$. Then,

\vspace{2mm}

\noindent \begin{multline*} \mbox{The Real Gromov-Witten invariant
on} \hspace{1mm} \mathcal{C} \\:= \sharp \{p \mid p \in ev_{1}^{-1}(
\Lambda_{1}) \cap \ldots \cap ev_{k}^{-1}(\Lambda_{k}) \cap
\md{k}{\CP^2}^{re} \}
\end{multline*}
\end{Def}
\vspace{2mm}

\par In the $I_{d}$([point], \ldots, [point],[line], \ldots, [line])
case, gains or losses of transverse intersection points in
$ev_{1}^{-1}( \Lambda_{1}) \cap \ldots \cap ev_{k}^{-1}(\Lambda_{k})
\cap \md{k}{\CP^2}^{re}$ are from the two factors:
\begin{itemize}
\item nature of singularities in stable maps: See Theorem
\ref{t;cusp}
\item tangency conditions: marked points deformation in the
neighborhood of the stable maps whose marked points go to the
tangential intersection point with the complex lines
\end{itemize}

In the second case, the tangential order equals to the intersection
multiplicity at that point. The non-transversal property caused by
the tangency condition prevents the Welschinger's invariant from
being extended to the
$I_{d}$([point], \ldots, [point],[line], \ldots, [line]) case.\\

\par Thorough characterization of walls and chambers is possible in
the $I_{d}$([point], \ldots, [point],[line]) case. The
non-invariance of the Welschinger type invariant becomes obvious.

\begin{example}$I_{d}$([point], \ldots, [point],[line]) case:
\par Let's consider $\CP(d_{1}, \ldots, d_{3d})^{re}$, where
$d_{i}=0$ if $i \leq 3d-1$ and $d_{3d}=1$. Theorem \ref{t;cusp}
shows that $\pi^{-1}(\mbox{wall})$ constitutes a part of walls in
$\CP(d_{1}, \ldots, d_{3d})^{re}$, where $\pi: \CP(d_{1}, \ldots,
d_{3d})^{re} \rightarrow \CP(d_{1}, \ldots, d_{3d-1})^{re}$.
\par Let $l$ be a line in $\CP^2$. Let's consider the following subset
 of $\md{0}{\CP^2}$:

\begin{multline*}
  P_{l} := \{ [(f, \CP^1, a_{1}, \ldots, a_{3d-1})] \in \md{0}{\CP^2}
  \mid\\
 \sharp (f(\CP^1) \cap l) = d-1 \hspace{2mm}
 \mbox{and $f$ is an immersion} \}\end{multline*}

The image of the stable map represented by each element in $P_{l}$
meets $l$ tangentially with an intersection multiplicity 2 at one
point and transversally at all other points. Proposition
\ref{p;dense} and ~\cite[p114, Lemma(3.45)]{ham} imply that $P_{l}$
is a codimension one subvariety in $\md{0}{\CP^2}$. Let $F:
\md{3d-1}{\CP^2} \rightarrow \md{0}{\CP^2}$ be a forgetful map and
$ev := ev_{1} \times \ldots \times ev_{3d-1}: \md{3d-1}{\CP^2}
\rightarrow \CP^2 \times \ldots \times \CP^2 := \CP(d_{1}, \ldots,
d_{3d-1})$ be a product of $i$-th evaluation maps. Then, one can
easily see that $\mathcal{WNS}_{tan:l}:= ev \circ F^{-1}(P_{l})$ is
a codimension one locus in $\CP(d_{1}, \ldots, d_{3d-1})$ because
$ev$ is a local isomorphism along $F^{-1}(P_{l})$ and $F$ is a
submersion. Let's consider the following subset $\mathcal{WT}$:
\begin{multline*}
\mathcal{WT} :=\{\mathcal{WNS}_{tan:l} \times [l] \subset \CP(d_{1},
\ldots, d_{3d}) \mid [l] \in \CP(d_{3d}) \}
\end{multline*}

$\mathcal{WT}$ may be considered as a codimension one fibration in
the trivial fibration  $\CP(d_{1}, \ldots, d_{3d})$ over
$\CP(d_{3d})$. Obviously, the real part of $\mathcal{WT}$ is
non-empty. And the real part $\mathcal{WT}^{re}$ constitutes the
rest of the walls in $\CP(d_{1}, \ldots, d_{3d})^{re}$.
\end{example}

\section{Transversality properties on $\md{3d-1}{r \CP^2}$ and $\overline{M}_{k}(\CP^1 \times \CP^1,
(a,b))$,
and their real enumerative implications} \label{s;rcp2} \vspace{2mm}

Let $r \CP^2$ be $\CP^2$ blown-up at $r$ points. $r \CP^2$  is a
non-convex variety. That is, $H^{1}(\CP^1, f^{*} T r \CP^2) \neq 0$
for some stable maps $f$ to $r \CP^2$. Nevertheless, if we consider
the divisor class of $d \cdot$[line], then the real Gromov-Witten
invariants can be defined as in the case of $\md{3d-1}{\CP^2}^{re}$.
\\

\begin{lemma}\label{l;last}
Let $p: r \CP^2 \rightarrow \CP^2$ be a natural blow-down map.
Consider a morphism \textbf{P}$: M_{0}(r \CP^2, d) \rightarrow
M_{0}(\CP^2, d)$ defined by $[(f, \CP^1)] \mapsto [( p \circ f,
\CP^1)]$. Then, \textbf{P} is an embedding.
\end{lemma}

Proof. Consider the short exact sequence of sheaves:

\begin{center}
$0  \rightarrow T \CP^{1} \rightarrow f^{*} T r \CP^2 \rightarrow N
r \CP^2 \rightarrow 0$
\end{center}

The same calculation we have done in Lemma \ref{l;tan} shows that
the tangent space at $[(f, \CP^1)]$ is isomorphic to $H^{0}(\CP^1, N
r \CP^2)$. Since $f(\CP^1)$ represents the divisor class $d
\cdot$[line], $f(\CP^1)$ never intersects with the exceptional
divisors in $r \CP^2$. Thus, the degree of the locally free
(sub)sheaf of $N r \CP^2$ and $N$ is the same, where $0 \rightarrow
T \CP^1 \rightarrow f^{*} T \CP^2 \rightarrow N \rightarrow 0$. It
implies that \begin{multline*}
 T_{[(f, \CP^1)]}\md{0}{r \CP^2}
\simeq H^{0}(\CP^1, N r \CP^2) \\
\simeq H^{0}(\CP^1, N) \simeq
T_{[(p \circ f, \CP^1)]} \md{0}{\CP^2}
\end{multline*}

The Lemma follows. \hfill q.e.d.\\

\par The morphism in Lemma \ref{l;last} doesn't change the type of singularities in the image
curves $f(\CP^1)$ and $p \circ f(\CP^1)$. Thus, the codimensions of
the tacnode stable maps locus, triple node stable maps locus and
cuspidal stable maps locus in $M_{0}(r\CP^2, d)$ are one. Lemma
\ref{l;last} implies that the index of the morphism $ev \circ
\textbf{P}$ is the same as the index of the $ev$ map defined on the
moduli space $\md{3d-1}{r \CP^2}$ before the compactification.
Therefore, the same intersection theoretic properties stated in
Lemma \ref{l;last}, Theorem \ref{p;gent} and Theorem \ref{t;cusp}
hold on $\md{3d-1}{r \CP^2}$. By ~\cite[Lemma 2.2]{gath}, the
Gromov-Witten invariant of $\md{3d-1}{r \CP^2}$ and that of
$\md{3d-1}{\CP^2}$ are the same. \cite[Lemma 4.3]{gap} shows that
the enumerative meaning of the Gromov-Witten invariant on
$\md{3d-1}{r \CP^2}$ comes from the number of nodal stable maps in
$M_{3d-1}(r \CP^2, d)$ which pass through the general points in $r
\CP^2 \setminus (\CP^{1}_{1} \cup \ldots \cup \CP^{1}_{r})$, where
$\CP^{1}_{i}$, $i = 1, \ldots, r$, are exceptional divisors. Walls
and chambers can be constructed on $M_{3d-1}(r \CP^2, d)$ and
$3d-1$-fold product of $r \CP^2 \setminus (\CP^{1}_{1} \cup \ldots
\cup \CP^{1}_{r})$ in the same way we did in section \ref{ss;point}.
Since $f(\CP^1)$ and $p \circ f(\CP^1)$ are the same type of stable
maps, i.e., $f$, $p \circ f$ have the same type singularities, the
Welschinger's
invariant becomes invariant.\\

\par $\CP^1 \times \CP^1$ is a convex variety. Let's assume $a \neq 0$ and $b \neq 0$.
All results in the case of $\overline{M}_{2(a+b)-1}(\CP^1 \times
\CP^1, a \times b)$ can be reproduced by repeating the same
arguments we have done in section \ref{sec;cp2}.

\bigskip
\noindent{\bf \em Acknowledgments.} Although I am the one who wrote
all parts and did all down-to-earth treatments by myself, the
originator of this paper is Gang Tian who wanted to write this paper
as his own work. He is a hidden coauthor who gave up his authorship
at the final stage of this paper after he played his role as a
coauthor until that point because he didn't have enough time in
reading all details for proceeding his current project. The way to
establish a real version of Gromov-Witten invariants is from Tian's
ideas. The non-transversality properties on the cuspidal stable maps
locus and the reducible stable maps locus were predicted by him. I
appreciate Bill Fulton for his suggestion to contact Gang Tian. This
work was done when I was a postdoc fellow at MSRI. I thank MSRI for
providing nice working
circumstances.\\

\vspace{1cm}

\noindent Department of Mathematics\\
University of Montana-Western \\
710 South Atlantic Street\\
Dillon, MT 59725, USA\\

\noindent e-mail: $s_{-}$kwon@hotdawg.umwestern.edu

\end{document}